\documentclass[12pt, draft]{amsart}  
\usepackage[a4paper]{geometry} 
\usepackage{hyperref}
\usepackage{ifdraft} 
\usepackage{color}
\usepackage{soul}
\newtheorem{thm}{Theorem}[section]
\newtheorem{cor}[thm]{Corollary}
\newtheorem{lem}[thm]{Lemma}

\theoremstyle{definition}
\newtheorem{para}[thm]{\smallskip}
\newtheorem{defn}[thm]{Definition}
\newtheorem{notation}[thm]{Notation}
\newtheorem{rem}[thm]{Remark}

\newtheorem{question}[thm]{Question}

\newtheorem{problems}[thm]{Problems}

\DeclareMathOperator\Aut{Aut}
\DeclareMathOperator\GL{GL}
\DeclareMathOperator\SL{SL}
\DeclareMathOperator\PSp{PSp}
\DeclareMathOperator\Sp{Sp}
\DeclareMathOperator\PSL{PSL}
\let\union \cup
\let\intersect\cap
\let\includedin \subseteq
\newcommand{\call}{{\mathcal L}}
\newcommand{\nat}{{\mathbb N}}
\newcommand{\inv}{^{-1}}
\newcommand{\ffi}{{\mathbb F}}
\newcommand{\id}{\mathrm{id}}
\newcommand{\tp}{\mathrm{tp}}
\let\iso\simeq
\newcommand{\qftp}{\mathrm{qftp}}
\newcommand{\cale}{{\mathcal E}}
\newcommand{\op}{^{\mathrm{op}}}
\newcommand{\gen}[1]{\langle{#1}\rangle}
\newcommand{\diag}{\mathrm{Diag}}

\newcommand{\Gg}{\mathbb{G}}
\newcommand{\Gdag}{G^\dagger}

\newcommand{\calm}{{\mathcal M}}

\definecolor{orange}{rgb}{1,0.5,0}

\DeclareMathOperator{\Gtwo}{G_2}

\begin{document}
\date{\today}

\title{Pairs of separably closed fields and exotic groups}
\author{Zo\'e Chatzidakis }
\address{Universit\'e Paris-Cit\'e - Sorbonne Universit\'e, CNRS, IMJ-PRG}

 \author{Gregory Cherlin}
 \address{Rutgers University, Emeritus}
 \thanks{This work began at the Newton Institute in Spring 2005, in the
   context of a semester program on \emph{Model theory and applications
   to algebra and analysis.} Both authors heartily thank the Newton
   Institute for their support. Work of the second author supported in part by the National Science Foundation under Grant No. NSF-DMS-0100794.}

\begin{abstract}
We look at simple groups associated primarily with the general theory of Moufang buildings, and to analyze their relation to stability theory in the model theoretic sense. As it becomes quite technical in the details, a lengthy introduction surveys the developments at a less detailed level.

The text, beginning from the second section, first deals with some model theoretic algebra of fields, followed by an extended study of three associated families of simple groups coming from the theory of Tits buildings, Moufang polygons, and Timmesfeld’s theory of exotic analogs of $\SL_2$.

The field theoretic part is fundamental (\S\ 2). The rest of the paper relates this to group theoretic constructions, with two sections surveying the consequences for the original Tits and Timmesfeld theory before concentrating on the more exotic groups associated with Moufang polygons.

A good deal of the group theoretical material is expository, aimed to make the relevant structural information meaningful to those coming from the direction of model theory.

\end{abstract}

 \maketitle

\section{Introduction}
Our  aim here is to construct some simple stable groups which are not algebraic (hence, ``exotic'').  These are not, strictly speaking, ``new'' groups, but instances of a phenomenon discovered by Tits long ago, in connection with the classification of buildings of spherical type \cite{Tits}.
He called them groups of ``mixed type''.
We became aware of this much later, while looking into the classification of \emph{Moufang polygons} given in \cite{TW} and discussed below in \S\ \ref{Sec:Moufang}.
Moufang polygons can be classified broadly speaking into algebraic (associated with algebraic groups), classical (in a historical sense), and \emph{mixed}, reusing the term introduced by Tits to reflect both their similarities to the algebraic case, and the use of two fields rather than one in their construction; but in the case of Moufang polygons the meaning of the term becomes a bit broader.

So we have on the one hand the groups identified by Tits, which are analogs of algebraic groups in Lie rank at least $3$, but with a coordinatization involving two fields $k\le K$, and we have also various groups associated with Moufang polygons which are analogs of algebraic groups in Lie rank $2$, but associated with a considerably more intricate collection of coordinatizing structures (including some of Tits' original type, constructed from a pair of fields). There is also a rank $1$ analog of $\SL_2(K)$ due to Timmesfeld, which we will also consider. 

A very natural program is then the following:
\begin{enumerate}  
\item[(a)] Construct some stable algebraic structures of the sorts used by Tits, Tits/Weiss, or Timmesfeld.
\item[(b)] Deduce the existence of the corresponding stable simple groups. 
\end{enumerate}This turns out to be more subtle than appears at first.  So we aim not only to carry this through in some cases, but also to point out some issues that others might want to explore further.

In the Tits setting, things work out neatly but with more delicacy than one might expect. 
An ample supply of coordinatizing structures for Tits' purposes is afforded by Theorem \ref{thm1},
and in a generalized form, { by} Theorem~\ref{thm2}. We cover some cases relevant to the 
Timmesfeld construction and an interesting case from the Tits/Weiss classification. However one is not quite done at this point. 

One might expect that a general interpretability result would allow for the systematic treatment of 
step (b) above. This seems not to be the case
(see Question~\ref{qu:T0-def}). On the other hand, in the context of the groups of Tits' type associated to a pair of fields $k\le K$, this is the case. 

The problem in general is that when one moves beyond Tits' original setting, the groups are defined as those generated by a collection of subgroups. This is perhaps clearest in the rank 1 case (the Timmesfeld construction), which is given explicitly as a subgroup of $\SL_2(K)$ whose diagonal subgroup is generated by elements whose coordinates lie in an additive subgroup of $K$.
The situation in rank 2 is much the same, but the notation involved is a good deal more complicated.

In fact, one may take a slightly different point of view on all of this, one that emerges most clearly in the rank 2 setting (Tits/Weiss).  This becomes more technical.
We describe this now, but the reader might prefer to look first  at the more concrete rank 1 setting of \S\ \ref{Sec:Rank1} where everything can be worked out in detail, from first principles, and only then return to a consideration of rank $2$.

In any case, in the rank 2 setting, there are at least two groups naturally associated with a given Moufang polygon, and it becomes important to distinguish them, and to consider more generally the groups lying between them.
The first group is the full automorphism group of the Moufang polygon. 
The second group, called the \emph{little projective group}, is defined (by analogy with Chevalley groups) as the subgroup generated by the so-called \emph{root groups},  which are the fundamental building blocks of the group from the point of view of either the Chevalley theory or the theory of Moufang polygons, and in the classical cases are copies of the additive group of the field. 
These groups appear in the Moufang theory as subgroups of the automorphism group of the Moufang polygon, and then the group they generate is one of the main groups of interest within the automorphism group, and is certainly the smallest group of interest, for our purposes. 

In most cases the latter group is simple, and is the socle of the full
automorphism group (its unique minimal normal subgroup). Between this
group and the full automorphism group there are some other groups which are interpretable in the coordinate system for the group, and whose commutator subgroup is our simple group. So if we start with a stable coordinate system then we can associate a stable group with a simple socle to it, but in passing to the commutator subgroup,  while we gain simplicity, we may lose definability.

Accordingly, our exposition becomes more elaborate than we had expected, as we sort through these issues. To complicate matters, our sources for the three cases take varying points of view, from the explicit matrix theoretic point of view of Timmesfeld, to the style of Chevalley (and Steinberg) in terms of generators and relations in the Tits/Weiss setting, and (for the part that concerns us) much more directly in terms of the structure of algebraic groups in the Tits setting. So we have the choice of unifying our perspective or staying close to our sources as we go along. We try to unify the description, but at the same time we do need to quote specific material from each source. 

The paper is aimed at model theorists with an interest in a variety of related topics. We have arranged it as follows:

In \S\ \ref{Sec:SCF}, Theorems \ref{thm1} and \ref{thm2} give the supply of stable ``coordinate systems'' with which we work. This is self-contained and is closely related to well-known work on  the model theory of separably closed fields.
Here the first theorem serves as warm-up to the second, and provides enough information to deal with the groups of mixed type as originally considered by Tits.
We describe such groups in \S\ \ref{Sec:Mixed:Tits} and prove that we do
indeed get stable simple groups of this type by passing to the context
of separably closed fields and applying Theorem \ref{thm1}.

Now Theorem~{\ref{thm2}} is of interest because the algebraic
systems considered are the natural parametrizing systems for the groups
which interest us. At the same time, the groups themselves cannot always
be defined in a first order way from these structures. A point of
considerable technical interest is that in some cases, enriching the original parametrizing structures to richer structures of the same kind may make the group first order definable.

We then pass to the opposite extreme---rank $1$---in the following section, working out the details of Timmesfeld's construction and the consequences for the issues of definability and interpretability that concern us here. Everything can be seen very simply by repeating standard computations (either from the point of view of Chevalley theory, or from the point of view of elementary linear algebra in two dimensions).  The unsettling phenomenon of a conflict between the desired simplicity and the desired interpretability appears at this stage. One can say more precisely how the initial coordinate system should be expanded to make the simple group definable, but then the issue of stability has to be approached afresh, and the situation becomes much more complex. Perhaps someone will investigate this further.

The last three sections discuss the related groups of automorphisms of Moufang polygons at some length. At this point the notation becomes noticeably more burdensome.
Here we encounter everything that we have seen in the original Tits
construction together with the complications that became visible in rank
1---and not much else, fortunately, other than some rather specific
notation. At this point one needs to work rather concretely in the
notation of root systems in order to sort out the details.  Readers will
probably find our presentation either excessively terse or excessively
detailed, depending on their degree of familiarity with that
notation. The ultimate result, which is a theme throughout much of the latter part of Tits/Weiss---though not put in these terms---is that in rank $2$ one has to deal with two separate instances of the rank $1$ theory, and otherwise things are rather similar to the case of algebraic groups. 

In more detail, the content of the last three sections runs as follows:  
In \S\ \ref{Sec:Moufang} we give an overview of what is done by Tits and Weiss
in \cite{TW}, and the notation used. Their goal is to give a classification theorem in terms of concrete coordinatizations by algebraic systems. This background material discusses what is common across all cases prior to the introduction of coordinates. 

The next two sections then look into two particular cases of the classification of Moufang polygons as given by Tits and Weiss. The first concerns Moufang hexagons, where we encounter examples already noticed by Tits as rank 2 analogues of the algebraic group of exceptional type $\Gtwo$ (so, $\Gtwo(K, k)$ in his notation). The second, more subtle example,  treated in the last section, concerns the Moufang quadrangles of so-called ``indifferent type,'' which are those most closely related to the Timmesfeld construction in rank 1. Our summary of the situation, above, focuses on this case: this is the setting which inherits the specific difficulties associated with the rank 1 case.

The classification of Moufang polygons involves further families which could be investigated model theoretically; they tend to involve structure incompatible with stability, but compatible, in principle, with simplicity. The interested reader may want to look further in that direction, and in particular investigate the problem of building coordinate systems of the various types which are simple in the model theoretic sense.

We imagine that most readers will be interested either in looking into \S\ \ref{Sec:SCF} and taking
much of the rest on faith (particularly from \S\ \ref{Sec:Moufang} onward), or else taking \S\ \ref{Sec:SCF} on faith and looking into the following group theoretic issues (including the definability issues that arise). Either approach should be perfectly feasible. Most of what we have to say in the group theoretic part is intended to be expository, but it was not always evident where to find clear statements in the literature of the facts most directly relevant to the model theoretic issues.

Up to this point, we have been very vague about the details, in order not to become lost in them. In the remainder of this introduction we give a more precise account of  the main points (and the key definitions) concerning the original construction of Tits, the lower rank constructions of Timmesfeld and Tits/Weiss, and the role of the model theory of separably closed fields in the construction of stable coordinate systems of the appropriate types.

\vspace{-0.1 pt} 

\subsection{The Tits construction: $G(k, K)$ \cite[(10.3.2)]{Tits}}
\label{Sec:Intro:G(k, K)}

Tits constructs analogs of (abstractly) simple algebraic groups over algebraically closed fields, in certain very special cases, defined from a suitable pair of fields $(k, K)$ with $k\le K$.
The point of view taken is that of Chevalley, with a small twist. This relies on the description of these groups in terms of root systems and their Dynkin diagrams, which may be summarized very rapidly as follows. This is either a reminder, or a few points of reference for the discussion afterward.

We begin with the algebraic group $G(K)$, which in algebraic terms is a $K$-split simple algebraic group of adjoint type. The $1$-dimensional subgroups are isomorphic to the additive or multiplicative group of $K$. A maximal torus $T$ is a product of a certain number of copies of the multiplicative group of $K$; that number is the Lie rank. The copies of the additive group of $K$ invariant under  the action of $T$ are the \emph{root groups} (with respect to $T$); these are permuted by the group $W=N(T)/T$; the action of $W$ on the root groups can be identified with the action of a finite reflection group acting on real Euclidean space (a Coxeter group) and these are classified by the Dynkin diagrams of types $A$--$G$. The root groups then correspond to a finite set of vectors invariant under the action of $W$ (these vectors encode the homomorphisms from $T$ to $K^\times$ which gives the action of $T$ on the corresponding root group).

From the Dynkin diagram, or the root system and the action of $W$, one can recover the construction of the group from the field $K$; this is the description of $G(K)$ as a Chevalley group.
We will see this concretely in the case of rank 2 in 
\S\S\ \ref{Sec:Hexagon}, 
\ref{Sec:Indifferent}, where in the latter case the construction is a generalization of the one described by Tits, and additional complications arise.

For our purposes it is important that the roots will always have either one or two root lengths.
The setting for the Tits construction involves 
a simple split algebraic group of adjoint type over a field $k$ associated with a root system in which,
in fact, two root lengths occur. 
Furthermore we require the characteristic to be ``exceptional'' in a certain sense (in a familiar sense from the point of view of finite group theory, and explained by Tits in terms of s
\emph{special isogenies}, \cite[(5.7.3)]{Tits}). 
The restriction on root lengths means concretely that the Dynkin diagram is of type $B_n$, $C_n$, $F_4$, or $\Gtwo$, and the restriction on the characteristic then means that the characteristic is $2$ unless we have type $\Gtwo$, in which case the characteristic will  be  $3$.\footnote{Here the classification 
by Dynkin diagrams can  be treated simply  as a set of labels for the cases of interest, 
until we come down to the rank $2$ case. 
Tits mainly deals with the case of $F_4$ in \cite{Tits}; he is able to identify types $B_n$, $C_n$ with groups he has treated from another point of view \cite[(2), p. 204]{Tits}, and $\Gtwo$ is mentioned in passing but lies outside the scope of that monograph.}

In this setting, one fixes a second field $K$ with
\[K^p \le k \le K.\]

With $G(k)$ the original algebraic group, one builds a group $G(k, K)$ containing $G(k)$, and contained in $G(K)$, much as one might construct $G(K)$ as a Chevalley group.

Namely, we consider a Borel subgroup 
$B=TU$ with $T$ $k$-split, we extend  the groups $T(k)$ and $U(k)$ to groups $T(k, K)$ and $U(k, K)$ in a manner to be described momentarily, and then we set $N(k, K)=N(k)T(k, K)$, so that $N(k, K)/T(k, K)$ is isomorphic to the usual Weyl group $W=N(k)/T(k)$. The group $G(k, K)$ is then defined as the group generated by $B(k, K)$ and $N(k, K)$.

The group $U(k,K)$ is an exact analog of the maximal unipotent subgroup of a Borel subgroup from the point of view of Chevalley. Namely, $U(K)$ is generated by the root subgroups, which are copies of 
$K_+$, subject to the \emph{Chevalley commutator relations} determined by the root system. 
One adjusts this construction by taking  the root groups for long roots to be copies of the additive group of the smaller field $k$, and the root groups for the short roots to correspond to the larger field. One may then check that the Chevalley commutator formula makes sense (using the precise data in that formula, and the particular value of the characteristic).

At this point, one could reasonably proceed as follows: using the same modified notion of root group based on a pair of fields $(k, K)$, take the group inside $G(K)$ generated by all the long root groups over $k$ and the short root groups over $K$. However, Tits proceeds in a different way.
 which connects up directly with his theory of \emph{BN-pairs}.  Before following him on this path, we discuss why one might do that.
 
 \subsubsection{BN-pairs and the Bruhat decomposition}
 
 In the first place, Tits' BN-pair theory gives a direct route toward connecting the new groups with the subject of his monograph \cite{Tits}. In the second place, the data $B(k, K)$ and $N(k, K)$ are explicitly given direct analogs of the usual groups $B(K)$ and $N(K)$. On the other hand the group \emph{generated by them} is potentially obscure; a priori it might very well be $G(K)$, for example. But the BN-pair theory implies a so-called \emph{Bruhat decomposition} 
 \[ G=\bigsqcup_W BwB\]
which is the double coset decomposition of $G(k, K)$ with respect to $B(k, K)$. (More properly, $w$ is replaced by a representative in $G$, but the corresponding double coset is well-defined.)
Comparing this to the corresponding decomposition of $G(K)$, we see that $B(K)\cap G(k, K)$ is $B(k, K)$, which is reassuring. And more generally, the Bruhat decomposition can be read as saying that $G(k, K)$ is built from  $B(k, K)$ in exactly the way that $G(K)$ is built from $B(K)$.

\subsubsection{The groups}
The groups obtained in this manner are (in the Dynkin notation)  the families $B_n(k, K)$, $C_n(k, K)$, and the exceptional groups $F_4(k, K)$, $\Gtwo(k, K)$.  
The groups $C_2(k, K)$ are variations on the algebraic group $\PSp_4(K)$.
Further variations are possible: these correspond to Moufang quadrangles of indifferent type in the sense of Tits and Weiss, discussed in \S\ \ref{Sec:TWvsTim}.
Rather than taking a pair of fields $k, K$, we take a large field $K$ and two additive subgroups
$K_0$, $L_0$, with
\[K^2 \le L_0\le K_0\le K\]
where now $L_0$ is a vector space over $K^2$ and $K_0$ is a vector space over the field generated by $L_0$. We then proceed to build a group $\PSp_4(L_0, K_0)$ in the manner of Tits, using $L_0$ to parametrize long root groups, $K_0$ for short root groups. This is not the description used by Tits and Weiss however; they build its associated Moufang polygon and then compute the subgroup of the automorphism group generated by the corresponding root subgroups (parametrized by $L_0$ and
$K_0$ rather than $k$ and $K$).

\smallskip

We have some unfinished business to attend to. 
On the one hand, we need to complete the definition of the groups $G(k, K)$. 
On the other hand, we should say a bit more as to how one actually obtains the BN-pair properties, or at least the Bruhat decomposition; this is the only way one has of seeing that these groups are in fact new groups, and Tits refers a little vaguely to Chevalley for this point, in \cite{Tits}, though elsewhere he gave the argument explicitly (in the Chevalley context).

\subsubsection{$G(k, K)$ (definition, concluded)}

We have described $U(k, K)$ as the subgroup of $U(K)$ generated by modified root subgroups.
Tits defines the torus $T(k, K)$, as the subgroup of $T(K)$ whose elements act sensibly on the root groups: 
that is, the elements of $T(K)$ which leave the root groups of $U(k, K)$ invariant.
In other words, these are the elements which act
via multiplication by an element of $k$ on the long root groups.

In particular the group $T(k, K)$ normalizes the group $U(k, K)$,  
and so we can define a  ``Borel subgroup'' $B(k, K)=T(k, K)U(k, K)$. (This 
is the largest available torus inside $G(K)$.) One could consider other constructions defining the torus in a different way. 
In the rank $2$ case this point is the subject of extended calculations in \cite{TW}; however the full automorphism group also contains elements inducing automorphisms of the coordinate system, which in the cases of interest to us are certain field automorphisms, and these will not appear in an algebraic group.

It is then reasonably clear that the ``Borel subgroup'' $B(k, K)$ is
interpretable in the pair {$(K, k)$}; more concretely, its underlying
set is definable in $G(K)$ if we take $K$ to be equipped with a
predicate 
for the subfield $k$. It then follows from the Bruhat decomposition that the same applies to $G(k, K)$, and thus stability of the coordinate system will give rise to stability of the group; the converse also holds (indeed $(k, K)$ is interpretable in $U(k, K)$).

\subsubsection{$B$ and $N$}
We come back to the point that the groups $B(k, K)$ and $N(k, K)$ 
give a Bruhat decomposition for $G(k, K)$, indicating how this goes in the 
setting of Chevalley groups, and how it relates the theory of BN-pairs.
For brevity we will now write $G$, $B$, $N$, $U$, and $T$ for the various groups involved in the definition
of $G(k, K)$. So the Bruhat decomposition is
 \[ G=\bigsqcup_W BwB,\]
As Tits mentions in \cite[10.3.2]{Tits}, a key ingredient is the fact that the 
nilpotent group $U$
can be written as the product of its root subgroups, taken in any order. Another ingredient is the fact that the Bruhat decomposition holds in rank one (in $\SL_2(k)$, $\SL_2(K)$, or the projective versions of these groups). To this one adds some observations about the operation of the reflections corresponding to a simple root on the set o positive roots, and the fact that opposite root groups generate a rank one subgroup.

We run over some of the more formal aspects of this argument, taking as our initial goal the Bruhat decomposition. As $G$ is generated by $N$ and $B$, and $W=N/T$ with $T$ contained in $B$,
the double coset decomposition exhibited is contained in $G$ and in order to show that it is $G$, it suffices to show that it is closed under multiplication by (representatives for) $W$ and under multiplication by $B$, the latter point being evident.
Also, as $W$ is generated by reflections $w_\alpha$ corresponding  to simple roots $\alpha$,
it suffices to check that sets of the form $w_\alpha BwB$ are contained again in the double cosets exhibited. What is claimed, in fact, is the following:
\[w_\alpha BwB\includedin Bw_\alpha w B \cup BwB.\]
This is one of the fundamental axioms in the theory of BN-pairs, in fact, so the question is how 
to verify it. 

This can be further reduced by similar formal manipulations, 
since $B=TU$ and $W$ normalizes $T$,
to a consideration of $w_\alpha U w$, and then even further by consideration of the structure of $U$. Namely, $U$ may be written as $U^*U_\alpha$ where $U_\alpha$ is the root group corresponding to $\alpha$, and where $U^*$ is the product of the remaining root groups, which is itself invariant under $w_\alpha$. One reduces quickly to a consideration of $w_\alpha U_\alpha w$. 
Then either $w$ or $w_\alpha w$ carries $U_\alpha$ into another root group group contained in $U$. In the first case $w_\alpha U_\alpha w\includedin w_\alpha w U$ and one finds that 
$w_\alpha B w B= Bw_\alpha w B$. In the second case one applies the same reasoning to 
$w_\alpha w$ in place of $w$, but one also uses the Bruhat decomposition for the rank one group generated by $U_\alpha$ and $U_{-\alpha}$. 

The last details are found in  the proofs of
\cite[Lemma 25, \S\ 3; (b) p. 34]{St} or 
\cite[(16), p. 323]{Tits-BN}.

\subsection{Tits-Weiss and Timmesfeld: subtleties}
\label{Sec:TWvsTim}

So far, everything proceeds according to plan. Now complications arise as we encounter some variations corresponding to Lie ranks 1 or 2, where the underlying algebraic systems are of a more general type. 

For us, the most interesting case concerns  Moufang quadrangles of ``indifferent'' type, 
similar to the buildings associated with Tits' groups $C_2(k, K)$, but more general.
Most of the complexity of this case, as far as the model theory is concerned, 
can be traced back to the rank $1$ groups associated with simple roots in this setting, 
which turn out 
to be the groups 
 Timmesfeld calls $\SL_2(L_0)$ and $\SL_2(K_0)$ (we are in characteristic $2$, so we do not need to distinguish $\SL_2$ and $\PSL_2$). 

There are interesting comments about the history and the differing emphases
of the various 
approaches taken  to this subject by \cite{Tits}, \cite{TW}, and \cite{Ti} to be found in 
Richard Weiss' review of \cite{Ti} in the AMS Bulletin
\cite{WeissBAMS}. In particular, the following has considerable relevance here:
\begin{quotation}
In a spherical building, groups of rank one appear as groups generated by pairs of ``opposite'' root groups, \dots. 
In the classification of Moufang buildings, in fact, these subgroups are
avoided to the maximal extent possible. The philosophy of abstract root
groups is just the opposite---groups of rank one are enshrined in the
hypothesis themselves and play a central role in the whole theory.
\end{quotation}
We will approach the rank 2 case via the rank 1 case, in order to encounter the model theoretic issues in their simplest ``pure'' state.
This means in particular that we will be crossing over between two rather different points of view.
\smallskip

We are again in characteristic $2$ with an imperfect field $K$, and we begin in rank $1$.
In the Timmesfeld setting---or rather, the  special case of interest to
us here---we will have an \emph{additive subgroup}  $L$ of $K$
containing $K^2$ and invariant under multiplication by
$K^2$. Timmesfeld's description of his group involves generation by two
``root subgroups'' parametrized by {$L$}, but as we will check
later, we can give a description similar to the one given by Tits
above.

We begin with a \emph{single} root group $U(L)$ (where $L$ is not necessarily a subfield) 
which we may take to be the upper unitriangular matrices with coefficients in the additive group $L$. 
If we followed Tits' construction we would also define a torus
$T(L)$ at this point. 
In fact we will take the root group $U(L)$ and its opposite, and the group they
generate, and then \emph{compute} the torus {$T(L)$} generated as a subgroup of $T(K)$. 
This turns out to be parametrized by
the multiplicative subgroup of $K$ which is \emph{generated by} the nonzero elements of $L$. 
This is the point at which nondefinability enters into the picture.

On other hand, after this detour we could start afresh and define $T(L)$ as the particular group of
diagonal matrices just mentioned, 
then define $B(L)=T(L)U(L)$, and let $\SL_2(L)$ be the group generated 
by $B(L)$ and 
a suitable Weyl group element. The usual Weyl group element
\[\begin{pmatrix}0&1\\-1&0\end{pmatrix}\]
will do (and we can omit the minus sign, as the characteristic is $2$).
This preserves the connection with the Tits construction; but we will in fact take Timmesfeld's 
definition as our point of departure.

As there is only one pair of roots, the field $K$ does not play much of a role here, and it could be replaced by the subfield $k$ generated by $L$.

On the other hand, the torus $T(L)$ is not the strict analog of the one considered by Tits.
The direct analog of Tits' $T(k, K)$ in this context would be the subgroup of 
the diagonal group $T(K)$ which normalizes $U$.  But this is $T(K)$, since $L$ is a vector space over $K^2$. So that torus would depend on the choice of $K$.

Notice that it is the small torus $T(L)$ which is a maximal torus
in the simple group  $\SL_2(L)$.
But in general it is the larger torus $T(K)$, defined in the manner of Tits, 
which is definable from the coordinate system, so here we have a definable group $T(K)\SL_2(L)$ with simple socle, stable if $(K,L)$ is, and the commutator subgroup of this group 
is simple, but not necessarily stable. 

All of this can be checked by direct computations which we will make, and which are the usual computations made over a field in the context of Chevalley groups. In particular one verifies the Bruhat decomposition in this context, and that leads to a proof of
the BN-pair axioms also in rank $2$ (carried out in a different way in \cite{TW}).

Turning to this rank $2$ case, let us call the group associated by Tits and Weiss to the coordinate
system $(K; L_0, K_0)$
$G_0(L_0, K_0)$. 
Namely, one defines $U(L_0, K_0)$ by strict analogy with the case of Chevalley groups,
as in the algebraic group $\PSp_4(K)$, with $L_0$ and $K_0$ parametrizing the long and short
root groups respectively, and using the Chevalley commutator relation to define the group law.

In an algebraic group setting one may then take the opposite group and the group they generate; or in the setting of Moufang polygons one may define the corresponding Moufang polygon (with some effort) and then consider the group generated by root subgroups. From this point of view one also computes the torus (with considerable effort in this setting). 
This gives a simple group which is not necessarily first order definable, 
because the torus itself is not necessarily definable,
 and in fact rank one groups of type $\SL_2(L_0)$ and $\SL_2(K_0)$ are involved. 
 The analysis of Tits and Weiss determines both the minimal torus
 (splitting the normalizer of the group $U$ as $T\cdot U$ in the
 corresponding simple group) and the maximal torus (giving a similar
 splitting, but in the full automorphism group of the Moufang
 polygon)\footnote{Tits and Weiss give in \cite[\S\ 37]{TW} a complete description of the
 automorphism group of the polygon for the various Moufang examples, which involves an ``algebraic part''
 and a subgroup coming from automorphisms of the field $K$; here, by
 {\em full automorphism group} we will mean the ``algebraic part''.}.

The result is that inside the automorphism group of the Moufang polygon, and above the group
generated by root subgroups, 
we have a family of groups, corresponding to a family of ``tori'' (in a very broad sense, allowing actions by field automorphisms on the coordinate system).

The smallest of these groups is simple but not necessarily definable over the coordinate system (in the first order sense), while the largest is rather too large for any of our purposes; 
but in between one can find a definable group whose commutator subgroup is the corresponding simple group (i.e., the associated torus is abelian). Here, definability refers to definability in the structure $(K;L_0, K_0)$.

In particular when the coordinate system is stable, the closest we come, in general, to building a stable simple group is to build a stable group with simple commutator subgroup.

On the other hand,  as yet we have \emph{no negative results} 
in the more challenging cases.  In particular we do not know whether some of the simple groups which are not interpretable in the associated algebraic systems might themselves be stable, for other reasons.

This last is not intrinsically a group theoretic question, since
the simple group  of interest is definable from  a  coordinatizing structure expanding
$(L_0, K_0)$ by the torus $T$ of the group and its action on the root subgroups 
(and conversely, this structure can be recovered from the group, if one is careful about the formulation).
The torus can be made a little more concrete as it is a product of 1-dimensional tori which can be taken separately and come from rank 1 subgroups of Timmesfeld's type.

So stability of the simple group is equivalent to  stability
of the structure $(L_0, K_0)$ together with the two 1-dimensional tori associated with the rank 1 subgroups corresponding to simple roots, and their actions on all the root groups.
In this sense, one can set aside the simple group and work with an expanded language of fields instead.

\subsection{Some model theory of fields}

A few introductory remarks about the model theory of fields are also in order, just to set the scene properly.
From our perspective, what was intriguing was the central role of imperfect fields in all of these constructions, and the known fact that separably closed fields have stable theories. 
This is what suggested the current line of investigation, and, in particular, our interest in the case of Moufang quadrangles of indifferent type.

The question as to whether every stable field is in fact separably closed is of long standing (see for example \cite{KrP-SFW}). This question has been placed in a broader framework by Shelah and others, and occurs now in a number of formulations generally all going by the name of \emph{Shelah's conjecture} for (e.g.) dependent fields \cite{HaHJ-SDF}. This  broader question is being actively pursued at present and leads into very different issues outside stability theory.
But certainly in the present state of knowledge the only definite source
of constructions of stable simple groups in which fields can be
interpreted will pass through the theory of separably closed fields. If
one enlarges the scope to simple unstable theories, then some other
constructions from the theory of Moufang polygons  would come into play, involving automorphisms and various semilinear  or quadratic forms.

\medskip

We turn now to the details, beginning with the model theoretic algebra
that produces a good supply of stable structures suitable for use as
coordinatizing structures, {with three cases: Timmesfeld's rank one
  groups $\SL2(L)$, and the two families of rank two groups $\Gtwo(k,
  K)$ and
the indifferent type for $\PSp_4$.}.
With that in hand we will take up the three sorts of groups of interest, starting with 
{Tits' theory over pairs of fields}, where matters are simplest in some fundamental sense (though with the usual apparatus of algebraic groups, root systems, and also BN-pairs all in the mix).
Then we pass to the rank 1 case as a relatively transparent context where real problems of definability arise, before coming finally to the most interesting case, Moufang polygons of indifferent type, where the groups to be constructed are stable, with simple socle equal to the commutator subgroup, and nonalgebraic.

\subsection{Main results of the paper}

As explained before, our aim was to study from a model-theoretic point
of view examples  of ``exotic groups,'' preferably simple ones.
We concentrated on three cases:
$\SL_2(L)$, $\Gtwo(k, K)$ and  groups obtained as automorphism groups of the Moufang polygons
of \cite{TW} coordinatized by an \emph{indifferent set}, and in particular the groups generated in that setting by the \emph{root groups} associated with  the Moufang polygon (and a fixed apartment).

The first point is that stable coordinatizing systems exist in all three
cases. {This is proved by fixing an imperfect separably closed
  field $K$ of the
  appropriate characteristic, and studying their model theory in various
enrichments,  by subfields of $K$, or by $K^2$ vector spaces between $K^2$
and $K$. The results are valid in arbitrary characteristic, and the main
result in that section is:}

\medskip
\noindent\textbf{Theorem \ref{thm2}}. 
{\em 
Let $K$ be a separably closed field in characteristic $p$, and
     \[K^p=K_0\leq K_1\leq K_2\leq\cdots\leq  K_m\leq K_{m+1}= K\] 
  a chain of subfields of  $K$ containing $K^p$. Furthermore, 
   for $1\le i\le m$ let $R_i$ be an additive subgroup of $K_{i+1}$ which contains $K_i$ and is a   
  vector space over $K_i$, and which satisfies, in addition, the following two conditions:
 \begin{enumerate}
 \item  $K_i=\{a\in K\mid aR_i=R_i\}$, 
 \item Any subset of $R_i$ which is linearly independent over $K_i$ is $p$-independent over $K_i$.
 \end{enumerate}

  Then the structure $(K, K_1, \ldots, K_m, R_1, \ldots, R_m)$ is stable, and the complete theory is given by the properties stated together with simple numerical invariants: the dimensions of both $R_i$ and $K_{i+1}$ over $K_i$, as finite values or the formal symbol $\infty$.
}

\smallskip
One also obtains  a variation of this result by slightly modifying the vector
spaces $R_i$ (Theorem \ref{thm3}).

These results will be applied in characteristic $3$ to a pair of fields and in characteristic $2$ to two fields and two additive subgroups meeting the additional requirements.

{Let us first start with two results on the groups $G_0(k,K)$ (``\`a
  la Tits''). }

\smallskip\noindent
\textbf{Theorem \ref{Thm:simpleG0kK}}.
\emph{ Suppose that $G(k)$ is of adjoint type (centerless) and split over $k$  
Then for $K\ne \mathbb{F}_2, \mathbb{F}_3$, the group $G_0(k,K)$ is simple.
}

\medskip
For the Tits groups, stability of the group is equivalent to stability
of the coordinatizing pair of fields, and we have

\medskip\noindent
\textbf{Theorem \ref{Thm:G(k,K)}}.
{\em Suppose  $G(K)$ is simple of type of type {$B_n$}, $C_n$, $F_4$, or $\Gtwo$ and $(K, k)$ is a pair of fields
with
\[K^p\le k\le K\]
and $p$ the appropriate characteristic ($3$ for type $\Gtwo$, and $2$ otherwise).

Then the following hold:
\begin{enumerate}
\item If the pair of fields
  $(K, k)$  is a stable structure, then the groups $G_0(k,K)$  and $G(k,K)$  {are} stable.

\item If $K$ is separably closed then $G_0(k,K)$ {and $G(k,K)$ are} stable  groups.
\end{enumerate}
}

\smallskip\noindent
{(The converse of item (1) is proved separately for type $\Gtwo$ and
$\PSp_4$, see Theorem~\ref{bidef:G2UKk} for
  $\Gtwo$, 
  and Theorem~\ref{Thm:C2:U(k,K)} for $G=\PSp_4$.)}

In particular for the case of groups of type $\Gtwo$ we {achieved
  our goal, and obtain} a class of 
automorphism groups of Moufang hexagons {which are both stable and
  simple}.

\medskip
Coming now to the rank one case (Timmesfeld's exotic simple groups of type $\SL_2(L)$), we have the following standard facts:

\medskip
\noindent\textbf{Theorem \ref{Thm:BN}}. 
{\em 
Let $K$ be an imperfect field $K$ of characteristic $2$ 
and $L$
an \emph{additive subgroup} satisfying
\[K^2 \le L\le K,\]
where $L$ is a vector space over $K^2$.
Let $T(L)\leq \SL_2(K)$ be the subgroup of $\SL_2(K)$
with coordinates in the multiplicative  subgroup of $K$ generated by
$L^*$. Let  {$B=T(L)L$} and $N=T(L)\gen{w}$. 

Then we have the Bruhat decomposition
\[\SL_2(L)=B\cup BwB.\]
In particular, {$L$} is the group of upper unitriangular matrices in
$\SL_2(L)$, and $T(L)$ is the diagonal subgroup. 

Furthermore, $\SL_2(L)$ is simple.
}

\medskip
{The definability theoretic properties of the group $\SL_2(L)$ are more subtle and lead us to 
consider a slight generalization $T\SL_2(L)$ where $T$ is a subgroup of the diagonal matrices
in a larger group $\SL_2(K)$ over a field. 
We may take $T$ to contain the diagonal matrices of $\SL_2(L)$.}

\medskip\noindent
{
{\textbf{Corollary \ref{Cor:410}}}. 
{\em Given a (slightly generalized) Timmesfeld group $T\SL_2(L)$, with additive group $L$
and torus $T$, 
there is a structure $(\tilde K,L,\bar T)$ with
$\tilde K$ a field and $\bar T$ a subgroup of $\tilde K$ such that the following are equivalent:
\begin{enumerate}
\item The group $T\SL_2(L)$ is stable.
\item The structure $(\tilde K, L, \bar T)$ with the field structure on $\tilde K$ and the additive and multiplicative subgroups $L$ and $\bar T$ is stable.
\end{enumerate}
In particular, when $T$ is the subgroup of $\SL_2(L)$ consisting of diagonal matrices (and $T\SL_2(L)$ is $\SL_2(L)$), the corresponding group $\bar T$ may be taken to be the
subgroup  of $\tilde K^\times$ generated by the nonzero elements of $L$.}}

\medskip
This is the point at which one realizes that $\SL_2(L)$ is likely to be undefinable in first order terms relative to its natural coordinatization by $(\tilde K,L)$, and examples falling under 
Theorem~\ref{thm2} confirm this.

     \medskip
We now deal with our second example, associated to hexagonal systems of
type 1/F, and which turns out to coincide with $\Gtwo(k,K)=\Gtwo_0(k,K)$. 
     
\medskip\noindent
    \textbf{Theorem \ref{Thm:G2(k,K)}}.
    {\em 
Suppose  
$(K, k)$ is a pair of fields
with
\[K^3\le k\le K\]
Then the following hold:
\begin{enumerate}
\item The group $\Gtwo(k,K)$  is stable if and only if the pair of fields
$(K, k)$  is a stable structure.
\item If $K$ is separably closed then $\Gtwo(k,K)$ is a stable simple group.
\end{enumerate}
}

\medskip\noindent
    \textbf{Theorem \ref{bidef:G2UKk}}. {\em
 Let 
$(K, k)$ be a pair of fields in characteristic $3$ 
with
\[K^3\le k\le K\]
and let $U=U(k,K)$ in the sense of $\Gtwo(k,K)$.
Then each of $U$ and $(K, k)$ is definable in the other.
}

    \medskip This immediately gives\\[0.05in]
\textbf{Theorem \ref{thm:hexagon}}. {\em The group $\Gtwo(k, K)$ is stable
  (model-theoretically simple, NTP$_2$, NSOP$_1$,
  \dots) if and only if the pair of fields $(K, k)$ is stable
  (resp. model-theoretically simple, \dots).}
    
    \medskip
    Now we turn to our real interest: the rank two case, and specifically automorphism groups of certain Moufang hexagons (\S\ 6) and Moufang quadrangles (\S\ 4).

      \medskip\noindent
         \textbf{Theorem \ref{Thm:C2:U(k,K)}}. {\em 
Let $(K; L_0, K_0)$ be a weak indifferent set 
and let $U$ be the group $U(k,K)$ in the sense of $\PSp_4(k,K)$.
Then each of $U$ and $(L_0,K_0,+,*)$ is definable in the other, where
\[a*b=a^2b\]
on $K_0$.
}

            \medskip\noindent
                \textbf{Theorem \ref{Thm:Indifferent:Definability:U}}. {\em 
Let $(K; L_0, K_0)$ be a weak indifferent set, $T(K)$ a maximal torus of $\PSp_4(K)$,
and $T$ a subgroup of $T(K)$ normalizing the group
$\PSp_4(L_0,K_0)$
and containing $T(K)\intersect \PSp_4(L_0,K_0)$.

Let $\calm$ be the structure
\[(K_0;L_0,+,T,\mu)\]
consisting of the group $K_0$ with the subset $L_0$, the abstract group $T$ with its multiplication,
and the following additional structure:
\begin{enumerate}
\item  the map $\mu:K_0\times K_0\to K_0$ defined by $\mu(a,b)=a^2b$;
\item actions of $T$ on $K_0$ and on $L_0$ which
correspond to the actions of $T$ on two root subgroups $U_\alpha$, $U_\beta$
with $\alpha,\beta$ the two simple roots, where $\alpha$ is short and $\beta$ is long.
\end{enumerate}

Then the group $G=T\PSp_4(L_0,K_0)$ is interdefinable with $\calm$.

In particular, $G$ is stable if and only if $\calm$ is stable. 
}

\section{Stable pairs of fields and related structures}
\label{Sec:SCF}

\subsection*{Results}

For our applications, we need to work with \emph{pairs} of fields, 
or with more general structures (but again, in pairs). 
But what can be done with pairs of fields can also be done, in the same way, with more than two nested fields,  and with the 
more general coordinatizing systems called \emph{indifferent sets.}
Our first result in this line will be the following, which we will need
in characteristic $2$ and $3$, and with $m=1$, so that we have two
distinct fields at our disposal:

  \begin{thm} \label{thm1} 
  Let $K$ be a separably closed field of characteristic $p>0$, 
  and let 
  \[K^p=K_0\leq K_1\leq K_2\leq\cdots\leq  K_m\leq K_{m+1}= K\] 
  be a chain of subfields of  $K$ 
  containing $K^p$, viewed as a structure with predicates for the fields.
  Then the theory of this structure is stable.
  
  Furthermore, this theory is axiomatized by the stated  properties together with  a specification of
  the dimensions $[K_{i+1}:K_i]$ (as finite values or the formal symbol $\infty$).
 \end{thm}
 
 The method of proof will pass through an elimination of quantifiers in an appropriate language---the language customarily used for quantifier elimination in separably closed fields, reviewed below, together with the appropriate unary predicates.

 This result already supports the Tits constructions, including some in rank 2, 
 notably in the case of $\Gtwo$, which was first described in \cite[\S\ 10.3, p. 205 (Remark)]{Tits}.
 
 But as we have explained, we need a more varied supply of coordinatizing structures, involving some additive subgroups as well as subfields---in characteristic $2$. The following will be sufficient for our current purposes, though as previously discussed, the question of stability of the associated simple groups would require even more elaborate coordinatizing structures, at this greater level of generality.

The relevant value of $m$ in the next theorem will be $2$, as we will be working mainly with the two additive groups $R_1$ and $R_2$.

\begin{thm} 
\label{thm2}

Let $K$ be a separably closed field in characteristic $p$, and
     \[K^p=K_0\leq K_1\leq K_2\leq\cdots\leq  K_m\leq K_{m+1}= K\] 
  a chain of subfields of  $K$ containing $K^p$. 
  Furthermore, 
   for $1\le i\le m$ let $R_i$ be an additive subgroup of $K_{i+1}$ which contains $K_i$ and is a   
  vector space over $K_i$, and which satisfies, in addition, the following two conditions:
 \begin{enumerate}
 \item  $K_i=\{a\in K\mid aR_i=R_i\}$, 
 \item Any subset of $R_i$ which is linearly independent over $K_i$ is $p$-independent over $K_i$.
 \end{enumerate}

  Then the structure $(K, K_1, \ldots, K_m, R_1, \ldots, R_m)$ is stable, 
  and the complete theory is given by the properties stated, together with simple numerical invariants: the dimensions of both $R_i$ and $K_{i+1}$ over $K_i$, as finite values or the formal symbol $\infty$.

\end{thm}


\subsection*\textbf{Algebraic preliminaries}

\begin{defn} Let $F\supset E$ be fields of characteristic $p>0$.
\begin{enumerate}
\item A subset $B$ of $F$ is \emph{$p$-independent}
  in $F$ if
 $[F^p[C]:F^p]=p^{|C|}$  for every finite subset 
$C$ of $B$; otherwise, it is said to be
\emph{$p$-dependent.} A maximal $p$-independent subset $B$ of $F$ is called a
\emph{$p$-basis} of $F$, and  one then has $F^p[B]=F$.

\item A subset $B$ of $F$ is \emph{$p$-independent over
  $E$} in $F$ if  
  $[EF^p[C]:EF^p]=p^{|C|}$ whenever $C$ is a finite subset of $B$.
  Note that if $E\supset F^p$, we could
  equally say: $B$ is $p$-independent in $E^{1/p}$.

  \item The {\em degree of imperfection} of the field $E$ is $e\in
    \nat\cup \{\infty\}$ such that $[E:E^p]=p^e$. Equivalently, it is
    the cardinality of a $p$-basis if $E$ has a finite $p$-basis, and the symbol
    $\infty$ otherwise. 
\end{enumerate}
\end{defn}

\pagebreak[2]
\begin{notation} Let $K$ be a field of characteristic $p>0$.
  \begin{enumerate}
\item For each $n>0$, we fix an enumeration $m_{i, n}(x_1, \ldots, x_n)$, $0\leq
  i<p^n$, of the $p$-monomials in $x_1, \ldots, x_n$, i.e., of all
  monomials on $x_1, \ldots, x_n$ where the exponents are between $0$ and
  $p-1$. Without loss of generality, $m_{0, n}(x_1, \ldots, x_n)=1$ for
  each $n$.

\item The
\emph{$\lambda$-functions }$\lambda_{i, n}$ on $K$ are defined in the following way:

\item  $\lambda_{i, n}(a_1, \ldots, a_n; b)=0$ if $a_1, \ldots, a_n$ is
not $p$-independent in $K$, or if $b$ is not $p$-dependent on $a_1, \ldots, a_n$ in $K$;  else, 
\item
the values of  the
$\lambda_{i, n}$  are  uniquely defined by the condition
\[b=\sum_{i=0}^{p^n-1}\lambda_{i, n}(a_1, \ldots, a_n; b)^pm_{i, n}(a_1, \ldots, a_n).\]

\item
Let  $\call$ be the language of fields
$\{+,-,\cdot,{}\inv, 0, 1\}$,  and let the language 
$\call_\lambda$ be $\call \cup\{\lambda_{i, n}\mid
n\in\nat, 0\leq i<p^n\}$.
Observe that the inverse of the Frobenius map is
  $\call_\lambda$-quantifier-free definable 
on $K^p$: if $b\notin K^p$, then 
$\lambda_{0, 1}(b; x^p)=x$. 

\item Let $B$ be a $p$-independent subset  of $K$.  For each $n$ and $i<p^n$, we
  denote by  $\lambda^B_{i, n}:B^n\times
K\to K$  the corresponding 
restriction of $\lambda_{i, n}$. If $a\in K$, we will say that the
  $\lambda^B$-functions are {\em well-defined at $a$} 
  when $a\in
K^p[B]$. Similarly, the iterates of the $\lambda^B$ are said to  (all) be well-defined at $a$
if $a\in K^{p^n}[B]$ for all $n>0$.

\item Suppose we have a  nested sequence  of fields
\[K_1\leq \cdots\leq K_m\leq K_{m+1}=K.\]
We define the language
  \[\call^m=\call_\lambda\cup \{K_1,\ldots, K_m\}\cup \{\lambda^{K_j}_{i, n}
  \mid 
  \mbox{$n\in\nat$, $0\leq i<p^n$, $j=1,\ldots, m$}\},\]
  where the $K_j$ are unary predicates
  for the subfields $K_j$, 
  and the  function symbols
  $\lambda^{K_j}$
  are interpreted as the usual $\lambda_{i, n}$ functions on the field
  $K_j$, 
  and $0$ outside. If $B_j$ is a $p$-basis of $K_j$, then
  $\lambda^{K_j,B_j}$ denotes the $\lambda^{K_j}$-functions restricted
  to $B_j^?\times K_j$. 
\end{enumerate}
  \end{notation}

We now collect some useful results, mostly classical 
(and trivial if the degree of imperfection of $K$ is finite). 
We  will give most of the proofs,  though briefly. 
More detailed proofs can be found at various points in \cite{B}, \cite{D} or \cite{Sr86}.

\begin{rem}\label{rem1}  
\ 

\begin{enumerate}

\item 
Let $E$ be a subfield of $K$. Then the following are equivalent:
\begin{enumerate}
\item $K$ is a separable
extension of $E$ 
\item $E$ is closed under the
$\lambda$-functions of $K$
\item the elements of any (or, some) $p$-basis of $E$ stay
$p$-independent in $K$. 
\end{enumerate}
In this case, the $\lambda$-functions of $E$
and of $K$ agree on $E$.
\smallskip

\item Let $B\subset K$ be $p$-independent.  Assume that the iterates
  of the $\lambda^B$-functions are well-defined at the element $a$ of $K$, and let $A_0$
  denote the set of these iterates. Then $\ffi_p(B, A_0)$ is closed under the
$\lambda^B$-functions. 
Hence
 $\ffi_p(B, A_0)$
has $p$-basis $B$, 
$K$ is a separable extension of $\ffi_p(B, A_0)$, 
and  $\ffi_p(B, A_0)$ is closed under the $\lambda$-functions
  of $K$.    
\smallskip

\item Let $E$ be a subfield of $K$ closed under the
$\lambda$-functions of $K$. 
Assume that $B$ is a $p$-basis of $K$
  such that $E\cap B$ is a $p$-basis of $E$.   Let $C\subset K$ be closed under the
$\lambda^B$-functions. Then  $E(C)$ is closed
under the $\lambda$-functions of $K$. 

Note that in general it is not true that if $A_1$ and $A_2$ are
$\call_\lambda$-substructures of $K$, then so is  the field
$A_1A_2$. For example, take $a_1, a_2, a_3, a_4$ $p$-independent, and consider
$A_1=\ffi_p(a_1, a_2)$, $A_2=\ffi_p(a_3, a_1a_2+a_4^p)$.
\smallskip

  \item  Let $E$ be a subfield of $K$ closed under the
    $\lambda$-functions of $K$, and let $a\in K$. 
    If $A$ is the closure
    of $E(a)$ under the $\lambda$-functions of $K$, then $A$ is countably
    generated over $E$.
\smallskip

    \item
The $\lambda$-functions of $K$ extend uniquely to the separable closure
$K^s$ of $K$.
\smallskip

 \item   Suppose  the subfield $E$ of $K$ is an $\call^m$-substructure of $K$, and
let   $a\in K$. Then the $\call^m$-substructure $A$ of $K$ generated by
$E(a)$ is countably generated over $E$.
\smallskip

 \item 
 Let $K^p\leq L\leq K$.  Let $B_1$ be a $p$-basis of $L$ over
    $K^p$, and $B_2$ a $p$-basis of $K$ over $L$. 
    Then $B_1\cup B_2^p$
    is a $p$-basis of $L$ and 
   \[L=\bigcap_{n\in\nat}K^{p^n}[B_1,B_2^p].\]

\item Let $K$ and $L$ be separably closed fields,  and
   $E\leq K, L$  an $\call_\lambda$-substructure.
Suppose that $K$ and $L$ are both saturated of the same cardinality
$\kappa$ with 
 $\kappa> |E|+\aleph_0$.
  Let $B$ be a $p$-basis of $K$
  such that $E\cap B$ is a $p$-basis of $E$, and let $B'$ be a $p$-basis
  of $L$ containing $E\cap B$. If $f:B\setminus E\to B'\setminus E$ is a
  bijection, then $f\cup \id_E$ extends to an isomorphism $K\to
  L$.

\end{enumerate}
\end{rem}

\pagebreak[2]
\begin{proof} 
\ 

(1) See \cite{B} or a similar (general) text.
\smallskip

(2) If $c$ is a $p$-independent
$n$-tuple in $K$ and  $a, b$ are two elements of $K^p[c]$, then
$\lambda_{i, n}(c; a+b)$ and $\lambda_{i, n}(c; ab)$ belong to the ring
generated over $\ffi_p[c]$ by the elements
\[
\mbox{ 
 $\lambda_{i, n}(c; a)$, $\lambda_{i, n}(c; b)$
  for $0\leq
i<p^n$.}\] 
    Moreover, $a^{-1}=a^{-p}(a^{p-1})\in K^p[a]$; this gives the first
assertion, and the second follows by (1). \smallskip

(3) By (2), $E(C)$ is closed under the $\lambda^B$-functions, and
the result follows by (1).\smallskip

(4) Let $A$ be as above, and extend a $p$-basis of $E$ to a $p$-basis
$B$ of $K$. 
Let $A_0$  be the set of $\lambda^B$-iterates of $a$.  As this set of functions is countable,
the set $A_0$ is
countable, and involves only countably many elements of $B$. That is, there
is a 
countable subset $B_0$ of $B$ such that all iterates of the
$\lambda^{B_0}$-functions are well-defined at $a$. 

Now by (3), $E(B_0A_0)$ is closed under the $\lambda$-functions of
$K$, and contains $A$.
\smallskip

(5) We know that for each $m$, $a\in K[a^{p^m}]$, so that there are
polynomials $f_m\in K[X]$, depending only on the minimal polynomial of
$a$ over $K$, such that $a=f_m(a^{p^m})$ for all $m$. Given a $p$-basis $B$
of $K$ the  polynomials $f_m$ determine
uniquely the values of the iterates of the $\lambda^B_{i, n}(a)$.\smallskip 

(6) For each $j=1,\ldots, m+1$, select a $p$-basis $B_{j}$ of $K_{j}$
such that  $B_{j}\cap E_j$ is a $p$-basis of $E_j$. Then, using (4), we
build an increasing sequence $A_i$ of subfields of $K$, where $A_0$ is
the  $\lambda^{B_{m+1}}$-closure of $\{a\}$, and for $i>0$, with
$i\equiv j\; \mathrm{mod}(m+1)$, $A_i$ is the $\lambda^{K_j,B_j}$-closure of
$A_{i-1}\cap K_j$. Since each $A_i$ is countable by (4), so is
$A'=\bigcup_{n\in\omega}A_n$, and by (4), since $E$ and $A'$ are closed
under the functions $\lambda^{K_j,B_j}$, so is $E(A')$.
\smallskip

(7) As $B_2$ is a $p$-basis of $K$ over $L$, $B_2^p$ is a $p$-basis of
$K^p$ over $L^p$, and therefore $L=L^p[B_1,B_2^p]$ and $B_1,B_2^p$ is a
$p$-basis of $L$. 

Observe now that
\[L=\bigcap_n L^{p^n}[B_1,B_2^p]\leq \bigcap_n K^{p^n}[B_1,B_2^p]
\leq
K^p[B_1,B_2^p]=L.\]
\smallskip

(8) This is a straightforward back-and-forth argument, using the stability
of the theory of separably closed fields of infinite degree of
imperfection. 

By (5), $f\cup \id_E$ extends to some $\hat f$ defined on $E(B)^s$. 
Assume now that we have an isomorphism $g:E_1(B)^s\to E'_1(B')^s$
extending $\hat f$, with $E_1, E'_1$ $\call_\lambda$-substructures of $K, L$
respectively, such that  $|E_1|<\kappa$ and $E_1\cap B$ a $p$-basis of
$E_1$.

Let $a\in K\setminus E_1$.  By (4), the $\lambda^B$-closure $A$ of $E_1(a)$ is countably
generated over $E_1$, and adding countably many elements of $B$ if
necessary, we may assume that  $A \cap B$ is a $p$-basis of $A$; let
$A_1$ be countable, 
closed under the $\lambda^B$, and such that $A=E_1(A_1)$.

By saturation of $L$,
there is some $A_1'\in L$ which realizes $g(\tp(A_1/E_1(A_1\cap B)^s))$, 
and as $A_1$
is independent from $ B$ over $A_1\cap B$,  
is separable over
$E_1[A_1\cap B]$, and $\tp(A_1/E_1(A_1\cap B)^s)$ is
stationary, it follows that $A_1'$ realizes $g(\tp(A_1/E_1(B)^s))$. 
This
proves one direction, and the other is symmetric. 
\end{proof}

\subsection*\textbf{Proofs}

\smallskip\noindent
We now give the proofs of Theorems \ref{thm1} and \ref{thm2}. 
We restate  Theorem \ref{thm1} in a more explicit form as follows:

\medskip\noindent\textbf{Theorem \ref{thm1}.\ }{ \em 
  Let $K$ be a separably closed field of characteristic $p>0$, 
  and let 
  \[K^p=K_0\leq K_1\leq K_2\leq \cdots\leq K_m\leq K_{m+1}= K\] 
  be a chain of subfields of  $K$ 
  containing $K^p$, viewed as a structure with predicates for the fields.
  Then the theory of this structure is stable.
  
  Furthermore, this theory is axiomatized by the stated  properties together with  a specification of
  the dimensions $[K_{i+1}:K_i]$ (as finite values or the formal symbol $\infty$),
  and admits elimination of quantifiers in the associated language
  $\call^m$, with the  predicates $K_j$ 
  and the functions $\lambda_{i, n}^{K_j}$ interpreted naturally.}

\begin{proof}
Since all $K_i$ contain $K^p$ and are contained in $K$,
the sequence 
\[K_0\leq K_1\leq \cdots \leq K_{m+1}=K\] 
is a series of  purely
    inseparable extensions.

Let $T_K$ be theory stating that the sequence of fields has the stated properties, and which, in addition,
specifies the degrees 
$[K_{j+1}:K_j]$. 
We show first that this theory is complete and allows quantifier elimination.

     Let $\cale=(E, E_1,\ldots, E_m)$ be an {$\call^m$-}substructure
    of $(K, K_1,\ldots, K_m)$. 
    Then 
    $K^p\cap E=E^p\subset E_1$, 
    each
    extension $K/E$, $K_j/E_j$  is separable, and for $j>0$, $K_{j}$ and
    $E_{j+1}$ are linearly disjoint over $E_j$, since $E_{j+1}/E_j$ is
    purely inseparable.

  We may assume that $(K, K_1,\ldots, K_m)$ is saturated of cardinality
$\kappa>|E|$, and we fix 
another model 
$(L, L_1,\ldots, L_m)$  
of $T_K$ 
containing $\cale$
which is also saturated of cardinality $\kappa$.

Let $B_1$ be a $p$-basis of $E_1/E^p$, $B_2$ a $p$-basis of
$E_2/E_1$, \dots, $B_{m+1}$ a $p$-basis of $E/E_m$; extend $B_1$ to a
$p$-basis $C_1$ of $K_1/K^p$, $B_2$ to a $p$-basis $C_2$ of $K_2/K_1$
(this is possible because $E_1$ and $K^p$ are linearly disjoint over
$E^p$, so that $K^p\leq K^pE_1\leq K_1$ are purely inseparable
extensions), \dots,
$B_{m+1}$ to a $p$-basis $C_{m+1}$ of $K/K_m$ (again, use $K/K_mE$ purely
inseparable). 

Then $\tilde C:=C_1\cup
\cdots \cup C_{m+1}$
is a $p$-basis of $K$, and for each $i=1,\ldots, m$, $\tilde C_i:=\bigcup_{j=i+1}^{m+1}C_j^p\cup
\bigcup_{j=1}^i C_j$ is a $p$-basis of $K_i$. Do the same with $L, L_1,\ldots, L_m$ to obtain
corresponding $p$-bases $D_1,\ldots, D_{m+1}$, and observe that necessarily, either $|C_i|=|D_i|$
is finite, or $|C_i|=|D_i|=\kappa$, by saturation of $K$ and $L$.

Thus, if $f$ is a bijection between $\bigcup_{j=1}^{m+1}(C_j\setminus B_j)$
and $\bigcup_{j=1}^{m+1}(D_j\setminus B_j)$ 
which sends each 
$C_j\setminus B_j$ to $D_j\setminus B_j$, 
then $\id_E\cup f$ extends to an
isomorphism of fields $K\to L$, which is the identity on $E$, and sends
each $K_j$ to $L_j$: use (8) and (7) in Remark \ref{rem1}. This shows that $T_K$ is complete and eliminates
quantifiers in $\call^m$.
\smallskip

Now we show stability of $T_K$. 
Let $a\in K$, and let $A$ be the
$\call^m$-substructure of $K$ generated by $(Ea)$
We may assume the  $p$-bases $C_i$ are chosen to contain a $p$-basis
of $A_i$ over $A_{i-1}$ extending $B_i$. 
By Remark \ref{rem1}(6) (and (4)), there is a countable
$\call^m$-substructure $A'$ of $A$ containing $a$ and such that 
 $EA'=A$ and $C\cap A'$ is a $p$-basis of $A'$. 
By
elimination of quantifiers, $\tp_{T_K}(a/E)$ is entirely determined by
$\qftp_{\call^m}(A/E)$, and because $A=EA'$ and by Remark \ref{rem1}(2), there are at most
$|E|^{\aleph_0}$ such types.
Thus the theory is stable.
\end{proof}

The proof of Theorem \ref{thm2} is very similar. Again, we reformulate
it in more precise terms:

\medskip\noindent
\textbf{Theorem \ref{thm2}.\ }\ 
{\em Let $K$ be a separably closed field in characteristic $p$, and
     \[K^p=K_0\leq K_1\leq K_2\leq\cdots\leq  K_m\leq K_{m+1}= K\] 
  a chain of subfields of  $K$ containing $K^p$. Furthermore, 
   for $1\le i\le m$ let $R_i$ be an additive subgroup of $K_{i+1}$ which contains $K_i$ and is a   
  vector space over $K_i$, and which satisfies, in addition, the following two conditions:
 \begin{enumerate}
 \item  $K_i=\{a\in K\mid aR_i=R_i\}$, 
 \item Any subset of $R_i$ which is linearly independent over $K_i$ is $p$-independent over $K_i$.
 \end{enumerate}
Then the structure $(K, K_1, \ldots, K_m, R_1, \ldots, R_m)$ is
  stable, and the complete theory is given by the properties stated
  together with simple numerical invariants: the dimensions of both
  $R_i$ and $K_{i+1}$ over $K_i$, as finite values or the formal symbol
  $\infty$, and admits elimination of quantifiers in the associated
  language $\call^m_R$, with the predicates $K_i,R_i$ and the
                      functions $\lambda_{i, n}^{K_j}$ interpreted naturally.}

\begin{proof}
Let 
\[\cale=(E, E_1,\ldots, E_m, F_1,\ldots, F_m)\subset (K, K_1,\ldots, K_m, R_1,\ldots, R_m)\]
 be a substructure. 
 As in
the proof of Theorem 
\ref{thm1}, the sequence 
\[E^p\leq E_1\leq \cdots\leq E_m\leq E\] 
is
purely inseparable, and each $K_i/E_i$ is separable. 

As usual, we
suppose that
the $\call^m_R$-structure 
$K$ is saturated, of cardinality $\kappa >|E|+\aleph_0$, and that
$(L, L_1,\ldots, L_m, S_1,\ldots, S_m)$ is another such model containing $\cale$.
By saturation, any of the invariants which are not finite take on the value $\kappa$ 
in both of the {$\call^m$}-structures $K$ and $L$.

The only change in what follows, relative to the proof of
Theorem \ref{thm1}, will lie in the initial choice of
$p$-bases $B_i$, $C_i$ and
$D_i$, so as to respect the additional structure.

Let $B_1$ be a $p$-basis of $E_1$ over
$E^p$ and extend it to a $p$-basis $C_1$ of $K_1$ over $K^p$.
For $i\geq 1$, 
let $B_i$ be a $p$-basis of $E_{i+1}$ over $E_i$, 
such that $\{1\}\cup (B_i\cap F_i)$ is an
$E_i$-basis of the $E_i$-vector space  $F_i$.
Extend $B_i$ to a $p$-basis $C_i$ of $K_{i+1}$ over $K_i$ in such
a way that $\{1\}\cup (C_i\cap R_i)$ 
is a $K_i$-basis of the
$K_i$-vector space $R_i$;
this is possible because 
$B_i\cap F_i=B_i\cap R_i$ is a $p$-basis of the purely
inseparable extension $E_i[F_i]$ of $E_i$, 
so that $K_i[F_i]\leq K_i[R_i]$ are also purely inseparable extensions of $K_i$.
Choose  $p$-bases $D_i$ within $L$ similarly.

As in Theorem \ref{thm1}, if $f_i:C_i\setminus B_i\to
D_i\setminus B_i$ is a bijection for $i=1,\ldots, m$,  then $\mathrm{id}_E\cup f_1\cup\cdots\cup f_m$ extends to an $\call_\lambda$-isomorphism $g:K\to L$,
which is an $\call^m$-isomorphism, and sends $R_i$ to $S_i$ for
$i=1,\ldots, m$. This gives completeness of the theory and also
  quantifier elimination for this language because of the
  $\lambda$-functions and conditions (1) and (2) on $(K_i,R_i)$.

The proof that the theory is stable goes much as before. 

Let $E$ and $K$
be as above, with $E=E^s$, and let $a\in K$. By Remark \ref{rem1}(4)(6), we know that there
is some countable $A_0\subset K$ containing $a$, 
closed under the
$\call^m_R$-functions, containing a $p$-basis of $A_0$, 
and such that $EA_0$ is closed under the
$\call^m_R$-functions. By stability of $(K, K_1,\ldots, K_m)$, there is
some countable substructure $E_0$ of $E$, which is separably closed, and
such that $\tp_{\call^m}(A_0/E)$ does not fork over $E_0$, and enlarging
$A_0$ we may assume that $A_0$ contains $E_0$ as an
$\call_R$-substructure. There are $2^{\aleph_0}$ 
possibilities for $\qftp_{\call^m_R}(A_0/E_0)$, and $|E|^{\aleph_0}$-many $\call^m$-formulas saying that
$\tp_{\call^m}(A_0/E)$ does not fork over $E_0$, so that there are at
most $|E|^{\aleph_0}$ types over $E$. 
Thus the theory is stable.
\end{proof}

As an easy corollary, we obtain

 \begin{thm}\label{thm3}
   Let $K^2=K_0\leq K_1\leq R_1\leq \cdots \leq K_m\leq 
   R_m\leq K$ satisfy the hypotheses of Theorem \ref{thm2}, and let
   $S_1,\ldots,S_m$ be additive subgroups of $K$, with $S_i$ a
   finite-dimensional $K_i$-vector space contained in $K_{i+1}$. Then the
   $\call^{2m}_R$-structure
   \[\mathcal{K'}=(K,K_1,K_1[R_1+S_1],\ldots,K_m,K_m[R_m+S_m],R_1+S_1,\ldots,R_m+S_m)\] is stable.

 \end{thm} 

 \begin{proof} As the $S_i$ are finite dimensional over $K_i$, both
   $R_i+S_i$ and $K_i[R_i+S_i]$ are definable (with parameters) in the
   $\call^{2m}_R$-structure $\mathcal{K}$.  \end{proof}

\begin{lem}
\label{Lem:DefinableFields}
Let $K$ be a field in characteristic $p$, and let
\[K^p=R_0 \leq  R_1\le R_2\le \cdots \le R_m \le K\] 
  an increasing chain of additive subgroups of  $K$.  
  Suppose that for all $i$ with $1\le i\le m$,
  $R_{i-1}\cdot R_i=R_i$. 
  Consider the structure 
  \[\calm=(R_m;R_0,R_1,\dots, R_{m-1},+, \mu)\]
  where the $R_i$ are given as subgroups of $R_m$ and 
  \[\mu:R_m\times R_m\to R_m\]
  is the function $\mu(a,b)=a^pb$.
  
  Then there are  \emph{$\calm$-definable }fields $K_0,K_1,\dots, K_m,\tilde K$ such that
 \[\tilde K^p=K_1\le R_1\le K_2\le \cdots \le K_m \le R_m\le K\le \tilde K\]
 and each $R_i$ is a vector space over $K_i$.
\end{lem}

\begin{proof}
The element $1$ in $R_m$ is clearly definable, hence the $p$-th power map $F:R_m\to R_m$
is definable. The restriction of multiplication to $R_m$ is a  partial binary operation $a\circ b$
defined by the relation
\[a*F(b)=F(c).\]
Define $K_m$ as the multiplicative stabilizer of $R_m$ under $\circ$:
\[\{a\in R_m\mid a\circ R_m=R_m\}
\]
This is a definable subfield of $R_m$ which contains $R_{m-1}$. 

Note that the structure on $\calm$ induces the corresponding structure on all $R_i$ and hence we have definable subfields $K_i\le R_i$ such that $R_i$ is a vector space over $K_i$, where in
addition $K_i$ contains $R_{i-1}$ if $i>1$, and $K_1$ contains $K^p$.
 Let $\tilde K= K_1^{1/p}$. Then $K\le \tilde K$. 
\end{proof}

From a model theoretic point of view, the reduced structure just on the $R_i$ is more 
convenient for interpretability results as the additional structure may then be treated as coming for free. Normally it would seem prudent, model theoretically, not to add undefinable structure to a given
coordinate system. In practice, that can be either highly undesirable or extremely convenient.
In the context of Theorem~\ref{thm2}, adding undefinable fields is harmless, and also at times extremely convenient. We will  see instances of the latter eventually (notably in the setting of rank 1 groups).

We conclude this section with some related questions.
which concern the choice of the vector spaces in Theorem \ref{thm2}, which gives a good  understanding of the most extreme case. Beyond that case, there may well  be other natural theories of similar kinds.

\begin{problems}
\ 

\begin{enumerate}

\item For $p>2$ and $K^p\leq K_1<K$, choose
  $(a_i)_{i\in\nat}$ $p$-independent elements of $K$ over $K_1$, and
  consider $R_1=\sum_{i\in\nat} K_1[a_i]$. Is $\mathrm{Th}(K, K_1, R_1)$
  stable? (The case $p=2$ is covered by Theorem \ref{thm2}). 

 Note that the union  $A=\bigcup_i \{K_1[a_i]\setminus K_1\mid i\in \nat\}$ is definable
 in the specified language via the formula
  $R_1(x)\land R_1(x^2)\land \neg K_1(x)$.  This tends to suggest a level of  complexity that may be incompatible with stability.
  In the structure as we have defined it, modulo $K_1$ all elements of $R_1$ 
  have a finite ``support'' in the set 
  $A$ (and in any model, elements with arbitrary finite supports will occur). 
  Whether this translates concretely into definable complexity remains unclear.

\item Let $K^p\leq K_1\leq K$ be separably closed fields of infinite
degree of imperfection, with $[K_1:K^p]=[K:K_1]=\infty$. 
Let
$\{a_i, b_i\mid i\in\nat\}$ be a subset of $K$ consisting of elements
$p$-independent over $K_1$, and set
\[R_1=\sum_{i\in\nat}K_1[a_i, b_i].\]
Is $\mathrm{Th}(K, K_1, R_1)$ stable?

 Note
that {again} a $\kappa$-saturated model will not be of the same form,
even if $p=2$. 

\end{enumerate}
\end{problems}

\section{Groups of mixed type $G(k, K)$ according to Tits \cite{Tits}; or a variation}
\label{Sec:Mixed:Tits}

In the present section we will discuss the groups of mixed type over pairs of fields
in the spirit of \cite{Tits} (with some slight variation). 
continuing on from the broad discussion in the introduction, 
\S\ \ref{Sec:Intro:G(k, K)}.
In this context, by applying Theorem~\ref{thm1}, we can identify some simple stable groups which are not algebraic but which one might reasonably call ``algebraic over two intimately connected fields.'' In Tits' monograph the focus was on rank at least $3$ as far as classification is concerned, but the constructions make sense in rank $2$, and in particular the case of $\Gtwo$ was covered
in 
\cite[\S\ 10.3, p. 205 (Remark)]{Tits}.

\begin{defn}\label{Def:G0kK}
Let $G(K)$ be a Chevalley group associated with a root system with roots of two lengths:
that is, type $B_n$, $C_n$, $F_4$, or $\Gtwo$. 

Fix a pair of fields $(k, K)$ satisfying
\begin{align*}K^p \le k \le K\end{align*}
where $p=3$ if the type is $\Gtwo$, and $p=2$ otherwise.

For $\alpha$ in the root system, define $U_\alpha(k,K)$ to be $U_\alpha(K)$ if $\alpha$ is short
and $U_\alpha(k)$ if $\alpha$ is long.

Let $G_0(k,K)$ be the group generated by the root subgroups $U_\alpha(k,K)$.
The $G_0$ notation indicates that we follow Tits' construction of
$G(k,K)$, but not exactly. The question is what part of the torus to take from $G(K)$ 
and as we will see in Lemma \ref{Lem:largetorus} that there is some latitude in this respect in the case of $C_2(k,K)$,
and more generally $\PSp_4(L_0,K_0)$.

\end{defn}

\begin{rem}
The group $G_0(k,K)$ has a BN-pair  \[B_0(k,K)=T_0(k,K)U(k,K), \ \ N_0(k,K),\]
where $U(k,K)$ is generated by root subgroups $U_\alpha(k,K)$ for
$\alpha$ positive, 
$T_0(k,K)$ is generated by the corresponding root tori, which can be defined as the intersection of the rank 1 group $\gen{U_\alpha(k,K),U_{-\alpha}(k,K)}$ with $T(K)$, or more directly as the groups $h_\alpha[U_\alpha(k,K)^*]$ in the notation of \cite[Lemma 19]{St}.
Then $N_0(k,K)$ may be defined as $N_{G(k)}(T(k))T_0(k,K)$ (which
normalizes $T_0(k,K)$ and has as quotient the Weyl group of
$G(K)$). That it constitutes a BN-pair can be proved with the
  classical arguments, using the fact that $U(k,K)$ can be written as a product of the root groups taken in any order,
and that the result holds for subgroups of the type of $\SL_2(k)$ and
$\SL_2(K)$ (treated more generally in \S\ \ref{Sec:Intro:G(k, K)}). \end{rem}

\begin{thm}\label{Thm:simpleG0kK}
Suppose that $G(k)$ is of adjoint type (centerless) and split over $k$  
Then for $K\ne \mathbb{F}_2, \mathbb{F}_3$, the group $G_0(k,K)$ is simple.
 \end{thm}

\begin{proof}
We use the Tits simplicity criterion for groups with a BN-pair, as
can be found in \S\ 29 of \cite{Hum}, see in particular Theorem~29.5. 

Since our groups have BN-pairs, 
it suffices to check  the  following points: 
\begin{enumerate}
\item [(a)]  $B$ is solvable and {centerless}.
\item [(b)] The set of generators of $W$ corresponding to the simple roots does not deompose into
a union of disjoint, nontrivial, commuting subsets.
\item [(c)] $B$ contains no nontrivial normal subgroup of the full group $G$.
\item [(d)]$G$ is perfect. 
\end{enumerate}

Of these four points, the first is the clear, and the second is a basic fact about 
the classification of the associated root systems. In terms of the usual Dynkin diagram representation
it means the diagram is connected. (In the rank two case with which we will be principally concerned,
it means that the two simple roots are nonorthogonal---so that the corresponding generators of the Weyl group do not commute.)

The third point may be argued as follows:
The group $B$ has a conjugate $B^w$ for which $B\cap B^w=T$, so any normal subgroup
$X$ of the full group
contained in $B$ would be contained in $T$. Then $[X,U]\le X\intersect U=1$ and $X$ centralizes $U$, forcing $X=1$ as the torus acts faithfully on $U$. This last point depends on the fact that the group has no center.

The proof that the group is perfect reduces to the condition 
$U_\alpha\le [U_\alpha,T]$ for the root subgroups $A$, since the root groups generate the full group.
This computation can take place in the rank $1$ group $\gen{U_{\alpha},U_{-\alpha}}$,
which is $\SL_2$ or $\PSL_2$ over one of the fields $k$ or $K$.
Here we may work concretely with $U_\alpha$ the group of strictly upper triangular matrices 
in $\SL_2$ and $T$ the group of diagonal matrices. 

Writing $x(a)$ for 
\[
\begin{pmatrix}
1&a\\0&1
\end{pmatrix}
\]
and $h(t)$ for the diagonal matrix with entries $(t,t^{-1})$, we have the commutator law
\[[h(t),x(a)]=x(a(1-t^{-2}))\]
Now we have only to choose $t$ so that $t$ is nonzero and $t^2\ne 1$ to get the general element
of $U_\alpha$ as a commutator.
\end{proof}

We gave the final computation explicitly as it will serve again in the more general setting
of Timmesfeld's rank one groups, below.

\begin{rem}\label{Rem:G2F3}
  There are no exceptions over $\ffi_3$, in fact, though for what one might call accidental reasons. Over $\ffi_3$, our definitions only allow one group, the algebraic group $\Gtwo(\ffi_3)$, and it is simple, for reasons like the ones we give but more delicate \cite[Lemma 32]{St}.

  \end{rem}

Of these, types $\Gtwo$ and $C_2$ 
recur below in the context
of Moufang  polygons (Moufang hexagons and Moufang quadrangles, respectively).
Type $C_2$ is a particular case of the class of Moufang quadrangles said to be of \emph{indifferent type}. As we will see, in a fairly precise sense, the class of groups associated to 
Moufang polygons of indifferent type is related to the narrower class of groups $C_{2,0}(k,K)$ in exactly the way that Timmesfeld's groups $\SL_2(L)$ are related to the usual groups $\SL_2(k)$ over fields.

From our point of view the interest of these groups lie in the
following:

\begin{thm}
\label{Thm:G(k,K)}
Suppose  $G(K)$ is simple of type of type {$B_n$}, $C_n$, $F_4$, or $\Gtwo$ and $(K, k)$ is a pair of fields
with
\[K^p\le k\le K\]
and $p$ the appropriate characteristic ($3$ for type $\Gtwo$, and $2$ otherwise).

Then the following hold:
\begin{enumerate}
\item If the pair of fields
  $(K, k)$  is a stable structure, then the groups $G_0(k,K)$  and $G(k,K)$  {are} stable.

\item If $K$ is separably closed then $G_0(k,K)$ {and $G(k,K)$ are} stable  groups.
\end{enumerate}
\end{thm}

We will look into this in a sharper form for the case of $\Gtwo$ in
\S\ \ref{Sec:Hexagon}. Clearly (1) relates to a couple of claims about interpretability and
(2) then follows via Theorem \ref{thm1}. But we give a proof of this form, in general.

It is also important to note that we set out with the expectation that something similar would occur 
in the analogous cases (at greater generality) in ranks 1 and 2, particularly in view of 
Theorem \ref{thm2}, but this is not the case: as already explained in the introduction, 
things become more subtle in rank 1 and then in rank 2 they remain equally subtle (but no worse).

\begin{proof}
In view of Theorem \ref{thm1} it suffices to prove the first point. For that we will use some coarse definability arguments. One should perhaps prove a bi-interpretability result characterizing definability exactly but it is not necessary for our purposes.

To show that the group is definable from the coordinate system (in first
order terms) we work inside the algebraic group $G(K)$, which is
certainly definable.  It suffices to show that the underlying sets of
{$G_0(k,K)$ and $G(k,K)$} are also definable, in the coordinate system $(K, k)$, as the group multiplication is
inherited.

In view of the Bruhat decompositions
\[G_0(k,K)=\bigsqcup_w B_0(k,K) w B_0(k,K)\ \ \hbox{and } G(k,K)=\bigsqcup_w B(k,K) w B(k,K)\]
with $w$ varying over a finite set of representatives, it suffices to show that 
$B_0(k,K)$ { and $B(k,K)$ are } definable.

Relative to the extended coordinate system $(K, k)$ the root groups are definable (parametrized by one of the fields).
The group $U(k,K)$ is the product (in any order) of its root subgroups, so it is definable.

The root tori that generate $T_0(k,K)$ are root tori of $G(K)$ or
$G(k)$, hence definable in the coordinate system. So $T_0(k,K)$ is
definable and $B_0(k,K)$ is definable. {The torus $T(k,K)$ is a definable
subgroup of $T(K)$ (in the pair of fields $(K,k)$), so $B(k,K)$ is definable.}\end{proof}

\begin{rem}
It turns out that the condition given in Theorem
  \ref{Thm:G(k,K)}(1) is also necessary: one  interprets the pair of
  fields $(K, k)$ in $G_0(k,K)$ and in $G(k,K)$, using the
  commutator relations. We will give the
  precise computations in two cases, see Theorem \ref{bidef:G2UKk} for
  $\Gtwo$, 
  and Theorem \ref{Thm:C2:U(k,K)} for $G=\PSp_4$. (In fact, in
  the case of $\PSp_4(k,K)$ we even prove outright definability.)
\end{rem}

Theorem \ref{Thm:G(k,K)} sets out the model for what we try to do in this paper.
This turns out to be more demanding than we initially expected. Theorem \ref{thm2}
prepares the ground by making an ample supply of  some coordinate systems needed to generalize Tits' construction in rank 2, but the definability issues are more severe as well. Namely, when one defines a group as ``the group generated by'' something,  and the coordinate system defines the generators, then the algebraist may be reasonably happy with that (particularly if a Bruhat decomposition results, and one can tell from that what group one has), but the model theorist needs to worry about the definability of the constituents of the Bruhat decomposition as well. One might reasonably object that if we had followed Tits we would also have defined not only the subgroup
$U(k,K)$ but the torus $T(k,K)$ and the group $N(k,K)$ as well, from the coordinate system, and the issue would disappear.  We will see  next why this is clearly not the case when we take up Timmesfeld's construction in rank 1, and then we will see why the difficulties that appear in rank 1 reappear in rank 2. 
{The only reason they do not appear in higher ranks is that
  the coordinate systems that appear in higher rank are of a
  particularly simple type, and in particular the only rank 1 groups
  that occur in that construction are $\SL_2(k)$, $\SL_2(K)$, and
  $\PSL_2(K)$.}

\section{The rank 1 case according to Timmesfeld \cite{Ti}}
\label{Sec:Rank1}

Timmesfeld presents a very general theory of groups generated by abstract root groups which includes the automorphism groups of most Moufang buildings, and starts off in rank 1 in what amounts to the study of split BN-pairs of rank 1 from another point of view. 
In particular, even the more exotic rank 1 groups arising as groups generated by pairs of opposite root groups in the context of Moufang buildings are captured by his theory.
We are interested in the ones which arise in the specific case of Moufang quadrangles of indifferent type, 
which we will come to in the next section.
In that case, we arrive at the particular rank 1 groups 
with which Timmesfeld begins his discussion in
\cite{Ti}, namely his Example~1.5,
as specialized further in \cite[Example~1.6 (2), p.~6]{Ti}.

In the presentation below, we begin with the explicit definition, but 
work out in detail the standard calculations in the manner of Chevalley or \cite{St}, in their minimalist form ($2\times 2$ matrices).
These calculations are identical to the usual calculations in $\SL_2(K)$, but we must pay close attention to where the entries of the matrices lie---and, in particular, which diagonal matrices are actually obtained in Timmesfeld's setting, and whether or not that set is  first order definable from the initial data.

\begin{defn}[$\SL_2(L)$ according to Timmesfeld]
We begin with an imperfect field $K$ of characteristic $2$ 
and an \emph{additive subgroup} $L$ satisfying
\[K^2 \le L\le K,\]
where $L$ is a vector space over $K^2$.
We then define the group $\SL_2(L)$ to be the subgroup of
${\mathrm SL}_2(K)$ 
generated by upper and lower unitriangular matrices in $\SL_2(K)$ 
with coefficients in $L$.

That is, we have the ``root groups'' 
$A$, 
$A\op$ consisting of the elementary matrices
\begin{align*}
a(t)&=\begin{pmatrix}
1&t\\ 0&1
\end{pmatrix}
&
b(t)&=\begin{pmatrix}
1&0 \\
t&1
\end{pmatrix},
\end{align*}
respectively, {with $t\in L$}. And we consider the group
$\SL_2(L)=\gen{A, A\op}$.
\end{defn}

There is a good deal to be said about the group $\SL_2(L)$. 
Our main concern is with a criterion for stability, which naturally leads us to consider related
definability issues, notable the definability of the subgroup of
diagonal matrices. This last issue turns out to recur substantially, afterward, in our discussion of rank 2 groups, as some of them contain 
Timmesfeld's groups. And for that matter, it is implicit in our treatment of Tits' construction
 (where we avoided beginning with a description of the torus), though in that construction the rank 1 tori involved were just the multiplicative groups of the two fields $k,K$. Here things become more delicate.

We begin with the Bruhat decomposition.
As a point of notation, 
we will denote by $L^*$ the set of non-zero elements of the additive group $L$. 
We make elementary calculations but keep track particularly of the diagonal matrices that appear.

\begin{thm} 
\label{Thm:BN}
Let $K$ be an imperfect field $K$ of characteristic $2$ 
and $L$
an \emph{additive subgroup} satisfying
\[K^2 \le L\le K,\]
where $L$ is a vector space over $K^2$.
Let $T(L)\leq \SL_2(K)$ be the subgroup of $\SL_2(K)$
with coordinates in the multiplicative  subgroup of $K$ generated by
$L^*$. Let  $B=T(L)A$ and $N=T(L)\gen{w}$. 

Then we have the Bruhat decomposition
\[\SL_2(L)=B\cup BwB.\]
In particular, $A$ is the group of upper unitriangular matrices in
$\SL_2(L)$, and $T(L)$ is the diagonal subgroup. 

Furthermore, $\SL_2(L)$ is simple.
\end{thm}

\begin{proof}  
Given any $a\neq 1$ in $A$, there is a unique
$b\in A\op $ such that $A^b=(A\op )^a$, and we  write $b=f(a)$; then  
$f(a(t))=b(-t\inv)$. (Even though we are in characteristic $2$, we use
 the minus sign
 since the computation works in any characteristic). 
 With $a_0=a(1)$, we find that
\[w:=a_0f(a_0)a_0=\begin{pmatrix}0&1\\
-1&0
\end{pmatrix}\]
is an element of the Weyl group of $\SL_2(K)$, and that the elements
$a(t)f(a(t))a(t)w$
are diagonal matrices in $\SL_2(L)$ of the form 
$\diag(t, t\inv)$, for $t\in L^*$. 

It 
follows that the subgroup of diagonal elements of $\SL_2(L)$ contains
all elements 
of the group $T(L)$.
From the formula 
$a(t)f(a(t))a(t)=\diag(t, t^{-1})w$, 
we deduce that \[b(-t\inv)=a(-t)\diag(t, t\inv)w a(-t),\]
so that 
$A\op \leq \gen{A, w,T(L)}$, and in fact,
 \begin{equation*}A\op \subseteq AT(L)wA \cup \{1\}.\end{equation*} 
 
Now we check that these calculations give
 \[\SL_2(L)=B\union BwB\]
 by formal manipulations, as in the case of fields.
 
 On the one hand, we know that both $T(L)$ and the element $w$ lie in 
 $\SL_2(L)$,  
 so the inclusion from right to left holds.
 In the opposite direction it suffices to check that the right hand side is closed under multiplication
 by $A$ and $A\op$, which is obvious for $A$. Hence for $A\op=wAw$ it suffices to check 
 closure under multiplication by $w$, which reduces to the following relations:
 \begin{align*}wBwB&=wT(L)AwB=T(L)wAwB=T(L)A\op B\\
 &\includedin T(L)(BwA\union \{1\})B
 \includedin 
 B\union BwB. \end{align*}
 
It now follows 
that $T(L)$ is the full diagonal subgroup of $\SL_2(L)$ and that $A$
is  
the full subgroup
of upper unitriangular matrices of $\SL_2(L)$, since this is clear in the case of the subgroup $B$,
and the double coset $BwB$ is disjoint from it.

For the simplicity of the group we use the BN-pair and follow the line of 
\cite[(16), p. 323]{Tits-BN}. 
We first show that the group $\SL_2(L)$ is perfect for $|L|>2$. It suffices to show that $A$ is contained
in the commutator subgroup, since then the conjugate $A\op$ is also contained in the commutator subgroup, and these two groups generate $\SL_2(L)$. 

We claim in fact that $[A,T(L)]=A$. We have 
\[[\diag(t^{-1}, t), a(s)]=a(s(1-t^2))\]
which for $t$ fixed and not equal to $0$ or $1$ represents a general element of $A$. The claim follows.

Now consider a normal subgroup $X$ of $\SL_2(L)$.

If $X$ is contained in $B$ then $X$ is contained in the conjugate $B^w$ and hence in the intersection, which is the group $T(L)$ of diagonal matrices in $\SL_2(L)$.
We then have $[X,A]\includedin X\intersect A=1$, so $X$ is in $C_{T(L)}(A)=1$, that is, $X$
is trivial.

So suppose now $X$ is not contained in $B$. Then the group $XB$ contains
$B$ properly, and is a union of $B$ double cosets, so by the Bruhat
decomposition $XB=\SL_2(L)$; 
 hence the quotient $\SL_2(L)/X$ is isomorphic to a quotient of $B$, and in particular is solvable. 
On the other hand as $\SL_2(L)$ is perfect the quotient is also perfect, and a perfect solvable group is trivial. So in this case $X=\SL_2(L)$. 

Thus $\SL_2(L)$ is simple.
For a statement from a broader point of view see \cite[I (2.10)]{Ti}.
\end{proof}

One should notice at this point that the torus $T(L)$ is likely to be undefinable in any natural language (at least, a priori; this is an interesting question in itself).
Accordingly, even if the structure $(K,L)$ is stable we run the risk that the group $\SL_2(L)$
is not. But there is a closely related group which is definable from the coordinate system,
and has $\SL_2(L)$ as its commutator subgroup: namely, the normalizer of $\SL_2(L)$ in $\SL_2(K)$. 
So we examine this.

\para\label{var-Tits}\textbf{The normalizer of $\SL_2(L)$}

For the present we fix the notation $K,L$ as in Timmesfeld's setting and consider
$\SL_2(L)$ within $\SL_2(K)$.

\begin{rem}
The full diagonal subgroup $T(K)$ of $\SL_2(K)$ normalizes $\SL_2(L)$, and the group
$T\SL_2(L)$ has the Bruhat decomposition
\[T\SL_2(L)=\hat B \union \hat B w\hat B\]
with $\hat B=T(K)A$.
\end{rem}

The point here is that diagonal matrices 
$\diag(t, t^{-1})$ 
act on $A$ and on $A\op$ by multiplication
by $t^{\pm 2}$, so $T(K)$ leaves $A$ and $A\op$ invariant. Then the Bruhat decomposition for 
$\SL_2(L)$ gives the Bruhat decomposition for $T(K)\SL_2(K)$.

The interest of this group is that it is definable over $(K,L)$ in view of the Bruhat decomposition,
and its commutator subgroup is $\SL_2(L)$ since $T(K)$ is abelian. Thus we have a definable stable group with simple commutator subgroup associated to any stable coordinate system $(L,K)$; this depends intrinsically on $K$ as well as $L$, though it would be very natural to take for $K$ the field generated by $L$ to get a more canonical construction (in similar settings in rank 2, this is actually part of the standard approach).

We note that in our definition of Tits' groups we preferred to follow Timmesfeld, and rather than defining a torus in advance, let it be computed in the group generated by root subgroups. As the coordinate system used was a pair of fields, the rank 1 subgroups $\SL_2(k)$ and $\SL_2(K)$ appearing there were not problematic. But we will need to keep these extra complications---and the need in some cases to sacrifice simplicity for definability---firmly in mind going forward.

\begin{lem}
In Timmesfeld's setting, the normalizer in $\SL_2(K)$ of $\SL_2(L)$ is $T(K)\SL_2(L)$.
\end{lem}

\begin{proof}
We work first in $\GL_2(K)$. Let $\hat T(K)$ denote the full subgroup of diagonal matrices.
This also normalizes $\SL_2(L)$. It suffices to check that the normalizer in $\GL_2(K)$
of $\SL_2(L)$ is $\hat T(K)\SL_2(L)$.
We have noticed that the normalizer contains this group.

Let $n$ belong to the normalizer of $\SL_2(L)$ in $\GL_2(K)$. 
If $n$ normalizes $A$ then it lies in the Borel subgroup $\hat T(K)A(K)$
of $\GL_2(K)$ (where $A(K)$ is the full set of  strictly upper
triangular matrices). Hence after multiplying by an element of $\hat
T(K)$ we may suppose $n\in A(K)$, and write $n=a(t)$ for $t\in K$.
In that case consider $L_1=\gen{L, t}$. Since $a(L_1)$ and $w$ lie in the normalizer of $\SL_2(L)$,
the group $\SL_2(L_1)$ is also contained in the normalizer of $\SL_2( L)$. But $\SL_2(L_1)$ is a simple group, so we find these two groups are equal and $n\in \SL_2(L)$.

If $n$ normalizes $A\op$ then $w n$ normalizes $A$ and we conclude similarly.

So suppose $A^n \ne A,A\op$. As the torus $\hat T(K)$
 acts transitively on the root groups of $\SL_2(K)$ (which correspond to the points of the projective line other than $0$, $\infty$) we may adjust by $\hat T(K)$ and suppose that $A$ is conjugated into a root group of the form $A(K)^b$ where $b\in \SL_2(L)$. 
 But then adjusting by this element of $\hat T(K)$ we may again take $n$ to normalize $A$,
and conclude as before.
\end{proof}

Thus the family of groups normalizing $\SL_2(L)$ in $\SL_2(K)$ is parametrized by the family
of groups $T_1$ lying between $T(L)$ and $T(K)$. We would like to take $T_1$ to be definable in 
$(K,L)$, ideally, but we would be perfectly happy as long as $(K,L,T_1)$  is stable.
Here $T_1$ is to be taken either as an abstract multiplicative group with an action on $L$ (corresponding to the action on $A$ in $\SL_2(L)$), or as the image of the action in $\Aut(L)$,
or more concretely as the multiplicative subgroup of $K$ whose action on $L$ is given by multiplication. Note that in the second interpretation the action of $\diag(t^{-1}, t)$
is multiplication by $t^2$  and in the third interpretation the multiplicative subgroup is actually 
the corresponding subgroup of $K^2$.

An attractive choice for the intermediate torus is the multiplicative group of the field
$K_L$ generated by
$L^*$. This will often not be definable in $(K,L)$, but we can work equally well with
$(K_L,L)$. And there are good chances that $T(L)$ will be equal to $T(K_L)$
in concrete cases; this leads to interesting questions.

Again: the choice of $T=T(L)$ gives a simple group; the choice of $T=T(K)$ gives a group definable in the original structure $(K,L)$ with $\SL_2(L)$ as commutator subgroup; and the choice
$T_1=K_L^\times$ gives a group which in general is not definable in $(K,L)$, but 
is definable in $(K_L,L)$; and if Theorem \ref{thm2} applies to $(K,L)$, it will also apply
to $(K_L,L)$. And as always, what we encounter here recurs in much the same form in rank 2.

We formalize the foregoing discussion further as follows:

\begin{thm} 
\label{Thm:Rank1:Definability}
Let $K$ be an imperfect field $K$ of characteristic $2$ 
and $L$
an \emph{additive subgroup} satisfying
\[K^2 \le L\le K,\]
where $L$ is a vector space over $K^2$.
Let $T$ be a group lying between the group $T(L)$  and the group $T(K)$.
Let $\bar T$ be
\[\{a\in K\mid \mbox{Multiplication by $a$ is induced by some element of $T$ acting on $A$}\}\]
Let $G=T\SL_2(L)$.

Then the following hold:
 \begin{enumerate}
\item 
\label{Thm:Rank1:Definability:G}
The group $G$ is definable in the
  structure $(L,\bar T,\cdot,\sigma)$, where $\cdot$ is the multiplication map on $L\times \bar T$,
  and $\sigma$ is the squaring map.
  
  \item Conversely, this structure is definable in $\SL_2(L)$.
  \end{enumerate}
  \end{thm}
  
  \begin{proof}
  \ 
  
  1. One builds the group $B$ from $\bar T$, $A$, and the action.
  One then builds the group $G$ as 
  \[B\union BwB=TA\union TAwA\]
  since $w$ normalizes $T$.  On the right side elements are uniquely represented either by
  pairs in $T\times A$ or by triples in $T\times A\times A$
   (since $A\op\intersect B=1$).
  Multiplication on this set is then determined by multiplication in $B$ and multiplication by $w$
  on the right. This is trivial for the map from $TA$ to $T A w$ and in the case of $TAwAw$ it reduces to the expression of $a(s)^w$ in terms of the Bruhat decomposition, given in the proof
  of Theorem \ref{Thm:BN} as
  \[b(-t)=a(-t^{-1})\diag(t^{-1}, t)w a(-t^{-1}).\]
  We may set aside the minus signs as superfluous. We need the operation of multiplicative inversion
  on $L$, which comes from squaring followed by the action of $\bar T$ on $L$, and the coordinate
  of $\diag(t^{-1}, t)$ in $\bar T$, which is $t^2$.

  2. $A$ and $T$ are, respectively, the centralizers  in $G$ of any of their nontrivial elements.
  So $G$ gives $A$ and $T$ and the action of $T$ on $A$. 
  This gives $\bar T$ as  a subset of $A$.
  
  The element $w$ allows us to define
  the function $f$ used in the proof of Theorem~\ref{Thm:BN}
   to compute the map
  from $a(t)$ in $A$ to $\diag(t, t^{-1})$ in $T$.  
  Thus we have the map from $a(t)$ to multiplication by $t^{-2}$ on $L$. This then gives both the set $\bar T$ as a subset of $L$, and its action on $L$ by multiplication. That is, $\bar T$ is the image of $a(1)$ under $T$, the image of $a(t)$ under the 
  corresponding element $\diag(t, t^{-1})$ of $T$ is 
  $a(t^{-1})$, and the squaring map is
 given by $a(t^{-1}\mapsto \diag(t^{-1}, t)\mapsto t^2$.
  \end{proof}
  
  \begin{cor}
  A group of the form $T\SL_2(L)$ in Timmesfeld's setting is stable if and only if the coordinatizing
  structure
  \[(L,\bar T,\cdot,\sigma)\]
  is stable.
  \end{cor}
  
  Now let us give a coordinatization that looks more normal from the algebraic point of view.
  
  \begin{thm}\label{thm:coord2}
  Let $K$ be an imperfect field of characteristic $2$, $L$ an additive subgroup of $K$,
  and $\bar T$ a multiplicative  subgroup of $K^2$ which contains $L^2$. 
  Suppose that ${\bar T}\cdot L\includedin L$.
  Then the structure
  \[(L,\bar T,\cdot,\sigma)\]
  in which $\cdot$ gives the multiplication on $\bar T$ and $\sigma$ gives the squaring
  map from $L$ to $\bar T$, is bi-interpretable with a structure
  \[(K_1,L,\bar T)\]
  where $K_1$ is a field satisfying Timmesfeld's conditions:
  \[K_1^2\le L\le {K_1}\]
  and $\bar T\includedin K_1^2$.
  \end{thm}
  
  \begin{proof}
  We have the multiplication on $\bar T$ and the squaring map to $T$.
  The restriction  $*$ of multiplication from $K$ to $L$ is given, as a partial function, by
  $a * b= c$ iff $a^2\cdot b^2=c^2$.
  
   Let $K_1$ be the multiplicative  stabilizer of $L$ in $L$:
   \[\{a\in K\mid aL\le L\}.\]
This is definable from $*$ and is a field containing $T$. Let $\tilde K$ be $K_1^{1/2}$ with its field structure, taken as an isomorphic copy of $K_1$ with an embedding $L\to \tilde K$ 
corresponding to the squaring map to $K_1$. We then have the structure
\[T\le \tilde K^2\le L\le \tilde K\]
with the multiplication on $\tilde K$ inducing the remaining structure.
\end{proof}

\begin{cor}\label{Cor:410}
In the (slightly generalized) Timmesfeld setting, the following are equivalent:
\begin{enumerate}
\item $T\SL_2(L)$ is stable.
\item The structure $(L,\bar T,\cdot, \sigma)$ with $\bar T\includedin L$, 
$\cdot$ the multiplication on $\bar T$, and $\sigma$ the squaring map from $L$ to $\bar T$, is stable.
\item The structure $(\tilde K, L, \bar T)$, with $\tilde K$ as above is stable.
\end{enumerate}
In the last clause, note also that the group
$\SL_2(L)$ in the sense of $K$ is also $\SL_2(L)$ in the sense of $\tilde K$) 
\end{cor}

One can do something quite similar in the rank $2$ indifferent case, in principle; namely there 
will be two rank 1 groups of Timmesfeld type and the condition is that both are stable (i.e.,
both exist within a \emph{single }stable structure).

Let us come back now to the case of $\SL_2(L)$, and consider the problem of stability.
This raises interesting questions of model theoretic algebra. We are considering structures
\[(K,L,T)\]
where $T=T(L)$ is the subgroup of $K^\times$ generated by $L^*$, given as an additional element of structure.  
By proper choice of $K$, the problem of stability for $\SL_2(L)$ becomes the problem of stability for structures of this kind. If $T$ is definable in $(K,L)$ there is no difficulty
(Theorem \ref{thm2}). If $T$ happens to be the multiplicative group of the field $K_L$ generated by $L$, then we may take the ambient field to be $K_L$, apply Theorem \ref{thm2} to that, and in this way force $T(L)$ to be definable. It is not yet clear how often that is the case. So the questions are of two sorts: when is $T(L)$ in fact the multiplicative group of a field, and in general, when is the expanded structure stable?

\begin{lem}
\label{Lem:Codim1}
Let $K$ be an imperfect field of characteristic $2$, 
and $L$ an additive subgroup of $K$ with
\[K^2\le L\le K\]
and $L$ a vector space over $K^2$. Suppose in addition that $L$ contains a subfield of codimension 
$1$ {(as a vector space over $K^2$)}. Then $T(L)$ is the multiplicative group of the field generated by $L$, and every element of 
$T(L)$ is the product of two elements of $L^*$.
\end{lem}

\begin{proof}
We write $L=K_1\oplus K^2 u$ for some $u\in L$, with $K_1$ a field.

Then $L$ generates the field $K_1(b)=K_1\oplus K_1b=K_1\cdot L\includedin L\cdot L$.
The claim follows.
\end{proof} 

Note therefore that we can always make $T(L)$ definable, in this setting by including the field $K_1$ in the coordinate system. In the context of Theorem \ref{thm2}, the theorem will continue to apply.

In particular, we have the following:

\begin{cor}\label{ex:T0-def}
 Let $K$ be an imperfect field of
  characteristic $2$, 
  with $[K:K^2]\geq 4$. Let $a, b$ be $2$-independent
  elements of $K$, and consider $L=K^2+aK^2+bK^2$. Then every nonzero element of
  $K^2[a, b]$ is the product of 2 elements of $L^*$, and therefore {$T(L)=K^2[a, b]^\times$}
  is definable in $(K,L)$. 
  \end{cor}

\begin{proof}
As $K_1=K^2(a)$ is a subfield of $L$ of codimension $1$, Lemma  \ref{Lem:Codim1} applies.
\end{proof}

Issues of stability in the groups $\SL_2(L)$  have led us to consider issues of definability in the underlying coordinate systems. It is clear that the field $K_L$ generated by $L$ plays a special role here. One has in general the question of definability of $K_L$ in some particular coordinate system, but when working in the context of separably closed fields, which is the only concrete case known currently, we have observed that this field should be added to the coordinate system and one should consider the issue of definability of $T(L)$ in the extended coordinate system, and in particular the question as to whether $T(L)$ always coincides with the multiplicative group of $K_L$, a question 
which reduces to the case of $T(L)$ finite dimensional over $K^2$. 

\begin{question}\label{qu:T0-def}  
Let $K$ be an imperfect field of
  characteristic $2$,  and $L$ an additive subgroup of $K$ containing $K^2$ which is
  a vector space over $K^2$.

\begin{enumerate}  
\item  Is $T(L)$ the multiplicative group of $K_L$? Does this hold at least if $K$ is separably closed?
\item If this is not the case, and the field $K$ is separably closed, is it possible for $T(L)$ to be definable in $(K, K_L)$ nonetheless?
\item Can $\SL_2(L)$ be stable when $(K,K_L)$ is not stable?
\end{enumerate}
  \end{question}

The first question, restricted to the case of $K$ separably closed, is the main question at present.
In the event of a negative solution, the second question should be taken as the natural refinement.
Finally, in a situation in which $T(L)$ is not definable in any structure covered by Theorem 
\ref{thm2}, the third question remains. This is not strictly a group theoretic question but a question about extending Theorem~\ref{thm2} to include certain multiplicative subgroups as well as additive subgroups, which seems very difficult. 

Since question (1)  reduces to the finite dimensional case and one can in principle make detailed computations in that case, it would be of interest to take up the minimal open cases, in which
$L$ has dimension $4$ over $K^2$, or more generally, where $L$ contains a subfield of codimension $2$. This seems accessible.

\section{The rank 2 case: Automorphism groups of Moufang polygons}
\label{Sec:Moufang}

Our own introduction to this subject came via the elegant work of Tits and Weiss in 
\cite{TW} concerning certain rank $2$ groups (or rather, the geometries on which they act).
So now we come, finally, to what was our point of departure. In practice we will focus on two of the cases which they consider, where the results of Theorem~\ref{thm2}, or the special case Theorem~\ref{thm1}, are directly applicable. In one case the group considered is the group $\Gtwo(k,K)$
already considered by Tits (though we give it a slightly different definition, one should bear in mind).
In the other case it is a substantial generalization of the Tits group of type $C_2$ in which the pair of fields used by Tits is replaced by a pair of suitably chosen abelian subgroups of fields in the manner of Timmesfeld.  

Here we run over the point of view of \cite{TW}, though as we find the groups easier to work with as subgroups of algebraic groups, we will adopt Tits' point of view for the more concrete discussions afterward. So this section indicates only how these groups were  originally identified, within the scope of a broad classification project (a project initially proposed in \cite{Tits} in a remark toward the end of the monograph).

The subject of \cite{Tits} is the theory of buildings, the geometries on
which simple algebraic groups, classical groups, and some other groups
act naturally; a classification is given in dimension at least $3$,
which can be taken as a classification of the corresponding
groups. These geometries generalize projective geometry, and just as high dimensional
projective geometries satisfy the Desargues condition and can then be
classified, all the higher dimensional buildings satisfy a related
Moufang condition, and are thus called Moufang buildings. Tits proposed
the problem of classifying all Moufang buildings in dimension $2$ or
higher; or more specifically, classifying them in dimension $2$
specifically and then  {reducing} the higher dimensional classification to that one. The project is carried through in \cite{TW}, with some surprises along the way.

In rank 2 the Moufang buildings are called \emph{Moufang polygons}. They are combinatorial
point-line geometries which are naturally represented as bipartite graphs where the parts are 
the points and lines, and the edge relation is incidence. One may also interpret the same graph with the points taken as lines and the lines taken as points, which would be treated as a dual geometry.
Accordingly the \emph{automorphisms} are taken to leave the points and lines invariant, and any graph automorphism which switches the parts would be called an \emph{anti-automorphism} in the geometric terminology. Tits and Weiss consider in great detail the structure of the geometric automorphism group $\Aut(\Gamma)$ of a Moufang polygon $\Gamma$ and in particular 
a certain subgroup $\Gdag$ which is almost always simple and which includes the usual Chevalley groups along with many other groups with a very similar structure. In particular, the theory begins with a definition of \emph{root subgroups} directly in terms of the action of the automorphism group on the graph, and $\Gdag$ is by definition the subgroup generated by a certain family of root groups (those associated with the vertices of an``apartment'', which is a cycle of minimal length in the graph). 

As in the case of Chevalley groups, one may define a ``maximal unipotent'' subgroup $U$
generated by half of the root groups (taking a path which covers half of
the cycle), which turns out to be a nilpotent group generated with the
root groups as generators and a generalized Chevalley commutator formula
as defining relations. 
We will consider some cases in which these commutator relations are the ones realized in some Chevalley groups.

Namely, we consider the Moufang hexagons which correspond to type $\Gtwo$ and more specifically to the groups $\Gtwo(k,K)$, and then the richer family of Moufang quadrangles of \emph{indifferent type}
 which correspond to type $C_2$, 
 and are realized in $\PSp_4$ (or $\Sp_4$ since we work in characteristic $2$). 
 In general, the polygon is called an $n$-gon if the shortest cycle length is $2n$: geometrically an $n$-gon has $n$ points and $n$ lines and forms a cycle of length $2n$ in the incidence graph. In particular the group $U$ is the (noncommuting) product of  $n$  root groups in a Moufang $n$-gon.

The main result of \cite{TW} is a classification theorem for Moufang polygons. Accordingly the various things known about Chevalley groups must not only be generalized, but proved in detail from first principles in a combinatorial setting. This complicates matters  relative to the theory
of Chevalley groups or algebraic groups, where the main facts are proved algebraically and may even be taken as belonging in part to the initial definition of the group (as in \cite{St}). 

But in addition to this, \cite{TW} contains detailed studies of the automorphism groups in all of the cases identified in the classification theorem, including that of the (mostly) simple group
$\Gdag$ as well as the full automorphism group and the quotient $\Aut(\Gamma)/\Gdag$, which can be viewed as a group of automorphisms of $U$. One of the main results of this analysis is the 
BN-pair structure for all of the groups between $\Gdag$ and $\Aut(\Gamma)$. As we have seen in the case of Timmesfeld's groups $\SL_2(L)$, we have reasons to consider larger groups than $\Gdag$ from the point of view of definability---though we will set aside the portion of $\Aut(\Gamma)$ which corresponds to nontrivial automorphisms of the coordinate system, which is not useful from the point of view of first order definability, and which does not appear in the corresponding algebraic group (when there is one).

\begin{rem}
\label{Rem:Uniqueness}
A very general lemma of \cite[(7.5)]{TW} states that a Moufang polygon is uniquely determined by the associated automorphism group $U$ and its sequence of root subgroups $U_1,\dots,U_n$. 
In particular the Chevalley commutator formula in $U$ determines the group $\Gdag$.
\end{rem}

We now describe the groups corresponding to the coordinate systems of indifferent type, which generalize Timmesfeld's systems $(K,L)$.

\begin{defn}
\label{defn:indifferent}
A \emph{weak indifferent set} is  a triple $(K, K_0, L_0)$, 
where $K$ is a field of characteristic $2$,
and $K_0, L_0$ are additive subgroups of $K$ for which
\begin{align*}
K^2 \le L_0\le K_0\le K,
\end{align*}
$L_0$ is  a vector space over $K_0^2$, and $K_0$ is a vector space over the field generated by $L_0$.

If a weak indifferent set satisfies the additional constraint that the field $K$ is generated by
the set $K_0$ then it is called an \emph{indifferent set}.
\end{defn}

It is customary to use indifferent sets in the strong sense in the literature, and we are introducing the terminology \emph{weak indifferent set} here to emphasize the variation. The distinction is not very significant from an algebraic perspective as there would be no harm in replacing the large field
$K$ in a weak indifferent set by the field generated by $K_0$. However, from a model theoretic point of view, the notion of weak indifferent set is axiomatizable, and the notion of indifferent set is not, so there is some advantage to allowing the broader notion into the formalism. It does not create any new examples of groups, however.

It is tempting to call a weak indifferent set an indifferent \emph{pair} (even though it is a triple)
because the groups $L_0,K_0$ play the roles previously played by the pair of fields $k,K$ in 
mixed type groups.
\begin{defn}
\label{Def:PSp4(L0,K0)}
Let $(K,K_0,L_0)$ be a weak indifferent set. Then
\[\PSp_4(L_0,K_0)\]
is the subgroup of $\PSp_4(K)$ generated by the subgroups
$U_\alpha(K_0)$ for $\alpha$ a short root, and by $U_\alpha(L_0)$ for $\alpha$ a long root.\
We call these groups the \emph{root subgroups} of $\PSp_4(L_0,K_0)$ (which will require a little justification).
\end{defn}

\begin{rem}
The group $\PSp_4(L_0,K_0)$ is defined by analogy with $\PSp_4(k,K)$, replacing
the pair $k,K$ by an indifferent set. As such it should more properly be denoted
\[\PSp_{4,0}(L_0,K_0)\]
and we may make use of that heavier notation if the point requires emphasis.

We have also identified a suitable torus in the verification of the BN-pair condition
and the Bruhat property, and so we could also have followed the route taken by Tits in defining
$G(k,K)$. But in any case it is important to us (and to \cite{TW}) that this group is generated by
its root subgroups.

There is some pathology in this construction, inherited from the rank 1
case, which will require close attention to the torus that appears in $\PSp_4(L_0,K_0)$, 
and to other tori that normalize this group.
\end{rem}

The definition of weak indifferent pair ensures that this group has more or less the same properties as
$G_0(k,K)$ where $G=\PSp_4$ and  $k=L_0$, $K=K_0$ are fields. We recall the relevant properties now.

First, the Chevalley commutator formula makes sense: that is, for positive roots $\alpha$, $\beta$,
and writing $U_\alpha$, $U_\beta$ for the root groups relative to $L_0$ or $K_0$ (as specified),
the formula giving coordinates of elements of 
$[U_\alpha,U_\beta]$ in the root groups of $\PSp_4(K)$ lie in the corresponding root groups
of $\PSp_4(L_0,K_0)$. However: this only works because in the special characteristics we consider, some terms in the general Chevalley commutator formula vanish, and the corresponding entries do not occur.  So this actually is what makes everything work.

At the same time, the  rank 1 groups $L_\alpha=\gen{U_\alpha,U_{-\alpha}}$ become
$\SL_2(L_0)$ or $\SL_2(K_0)$ in the sense of Timmesfeld.

One gets the BN-pair property, the corresponding Bruhat decomposition, and simplicity as previously. The computations we made in rank 1 close the gap between the usual $\SL_2(K)$
and the Timmesfeld variations, and the rest of the argument for the BN-pair is formal, modulo the rank 1 case.

Notice also that $\PSp_4(L_0,K_0)$ lies between $\PSp_4(K^2)$ and
$\PSp_4(K)$.

\begin{lem}\label{Lem:simpleG2Sp4}
 The groups $\Gtwo(k, K)$ and $\PSp_4(L_0, K_0)$ are simple (for $K$, $K_0
 \neq \ffi_2$).
\end{lem}

 \begin{proof}
 We use the Tits simplicity criterion for groups with a BN-pair, as
can be found in \S\ 29 of \cite{Hum}, see in particular Theorem~29.5.
Since our groups have BN-pairs, it suffices to check the following
points: (a) $B$ is solvable and centerless; (b) the set of generators of $W$ corresponding to the simple roots does not decompose into a union of disjoint, nontrivial, commuting subsets; (c) $B$ contains no nontrivial normal subgroup of the full group $G$; and (d) $G$ is perfect.

Of these four points, the first two are {clear} since there are only two simple roots and the corresponding reflections do not commute ($W$ is a dihedral group of order greater than $4$).
The other two points were noticed in the proof of the rank 1 case (Theorem~\ref{Thm:BN}), and the proofs given there continue to work. We repeat the main points.
The group $B$ has a conjugate $B^w$ for which $B \cap B^w= T$, so any normal subgroup of the full group contained in $B$ would be contained in $T$, after which it follows easily that it centralizes $U$, hence lies in $U$, hence is trivial.
The proof that the group is perfect reduces to the condition $A \leq [A, T ]$ for the root subgroups $A$, which is already shown in the rank $1$ case.\end{proof}

\begin{lem}
The groups $\Gtwo(k,K)$ and $\PSp_4(L_0,K_0)$ are the groups $\Gdag$ of
\cite{TW} corresponding to the Moufang hexagons of type $(1/F)$ and the Moufang quadrangles of indifferent type in the sense of \cite{TW}.
\end{lem}

\begin{proof}
We suppose the field $K\ne \mathbb{F}_2$.
By \cite[Thm.~6.1]{DemT}, if $G$ is the universal Steinberg group with the same presentation
as $\Gdag$ then $G/Z(G)$ is simple.

Since $\Gdag$ and the groups of type $\Gtwo$ or indifferent type are generated by root groups
satisfying the same relations, both are homomorphic images of $G$. Furthermore both groups
are simple  by Theorem \ref{Thm:simpleG0kK}, Remark \ref{Rem:G2F3}, Lemma \ref{lem:Gdag=G2kK} and \cite[(37.3)]{TW}. 
So the kernel in both cases is $Z(G)$ and the two quotients are isomorphic.
\end{proof}

\section{Some Moufang hexagons}
\label{Sec:Hexagon}

\para\textbf{The hexagonal case}

We return to a discussion of the groups of type $\Gtwo(k,K)$, which are the tamer examples of slightly exotic automorphism groups of Moufang polygons.  There are also the groups
$\PSp_4(k,K)$, which fall under the indifferent case treated in the next section, but we
will not single them out for attention.  So $\Gtwo(k,K)$ will serve as our model for the discussion when the definability issues important to model theory are not severe, and
the Theorem~\ref{thm1} (rather than the more general Theorem~\ref{thm2}) is adequate for our purposes.

Before going into details, we remind the reader that our definition of $\Gtwo(k,K)$ is not literally the same as Tits'.  On the other hand, as we are in characteristic $3$ it becomes quite clear that 
the torus Tits introduces in his definition is part of our group as well, and the two coincide.
In fact  there is only one possible torus in this case (as opposed to the situation encountered above in rank 1).

\begin{lem}\label{lem:Gdag=G2kK}
Working in $\Gtwo(K)$,
 the normalizer of the group $U(k,K)$ in the torus $T(K)$ is the torus
\[T(k,K) :=T(K)\intersect \Gtwo(k,K).\]
\end{lem}

\begin{proof}
The torus of $\Gtwo(k,K)$ contains the root tori $H_\alpha(k)$ for $\alpha$ long and $H_\alpha(K)$
for $\alpha$ short.
The torus $T(K)$ is the product of the two root tori for the simple roots $\alpha$ and
$\beta$ with $\alpha$ short and $\beta$ long.
So it suffices to check that the normalizer of $U_\beta$ in $H_\beta(K)$ is $H_\beta(k)$.

This is a rank 1 computation which would not work in characteristic $2$, but  we are
in characteristic $3$.  The action of the root group element $h_\alpha(t)$ (corresponding to
$\diag(t, t^{-1})$ in our rank 1 computations) is via multiplication by $t^{-2}$. 
We also have $K^3\includedin k$ so for any $t$ in $H_\beta(K)$ normalizing $U_\beta(k)$, both 
$t^{-2}$ and $t^3$ lie in $k$, and thus $t\in k$.
\end{proof}

One can make similar computations in other cases---in fact, whenever there are two long roots---but
working with rank 2 subgroups. On the other hand, this does not apply to $\PSp_4(k,K)$ and we will return to that point.

At this point we may use the notation $\Gtwo(k,K)$ with a clear conscience to refer to the group either as defined by Tits or as defined here.

Now we recall Theorem \ref{Thm:G(k,K)}, specialized to our context.

\begin{thm}
\label{Thm:G2(k,K)}
Suppose  
$(K, k)$ is a pair of fields
with
\[K^3\le k\le K\]
Then the following hold:
\begin{enumerate}
\item The group $\Gtwo(k,K)$  is stable if and only if the pair of fields
$(K, k)$  is a stable structure.
\item If $K$ is separably closed then $\Gtwo(k,K)$ is a stable simple group.
\end{enumerate}
\end{thm}

In one direction, the interpretability of the coordinate system in the group $G$ comes from
its interpretability in the Borel subgroup $B$; most of it would come in fact from the rank $1$
case where $B$ is just $K^\times$ acting on $K_+$, which is simple the field disguised as a group.
However when two fields $k,K$ are present one also has to embed $k$ into $K$ and for this one uses commutation relations in $U$. The question arises as to whether the commutation relations in $U$ are already enough to recover the coordinate system. This is obviously false in 
rank $1$ since $U$ is just an abelian group in that case. However it is to be expected (or hoped) in rank $2$, and we will make this analysis in the two cases of interest here.
In $\Gtwo$ life is greatly simplified (relative to the indifferent case) by the fact that  all root groups are parametrized by fields.

\begin{thm}\label{bidef:G2UKk}
Let 
$(K, k)$ be a pair of fields in characteristic $3$ 
with
\[K^3\le k\le K\]
and let $U=U(k,K)$ in the sense of $\Gtwo(k,K)$.
Then each of $U$ and $(K, k)$ is definable in the other.
\end{thm}

One major question is the extent to which the root subgroups of $U$ are definable.
As homomorphisms from $U$ to $Z(U)$ produce automorphisms of $U$ that move the root groups, 
this is not literally the case. Thus the absence of the torus is significant.

\begin{proof}
In one direction, in all of these constructions, the original definition of $U(k,K)$ is given in first order terms relative to the coordinate system, and there is really nothing to prove. We illustrate this in the case of $\Gtwo$ since in any case we require the definition of $U$ for the converse.

We have root groups $U_1,U_2,\dots,U_6$, with the roots alternately short and long, parametrized by $K$ or $k$
respectively. We form a group generated by the six root groups with the following commutator
relations:
\begin{align*}
\tag{[1,5]; [2,6]}[x_1(a),x_5(b)]&=x_3(-ab); \qquad [x_2(t),x_6(u)]=x_4(t u);\\
\tag{[1,6]}[x_1(a),x_6(t)]&=x_2(-t a^3)x_3(t a^2)x_4(t^2a^3)x_5(-t a),
\end{align*}
where 
$a, b\in K$ and $t, u\in k$, we also require that all other pairs $U_i,U_j$ commute. 
Recall that 
the characteristic is $3$ and some terms of the usual Chevalley commutator formula degenerate.

In particular the commutation formula shows that
\begin{align*}
Z(U)&=U_3\times U_4;&Z_2(U)=\gen{U_2,Z(U),U_5}.
\end{align*}

We set $\bar U=U/Z_2(U)$,  $\bar V=Z_2(U)/Z(U)$, and $Z=Z(U)$ and we have
\begin{align*}
Z&=U_3\times U_4;&
\bar V&= \bar U_2\times \bar U_5\iso U_2\times U_5;&
\bar U&=\bar U_1\times \bar U_6 \simeq U_1\times U_6.
\end{align*}

Furthermore commutation induces bilinear maps
\begin{align*}
\gamma&: \bar U\times \bar U \to \bar V&\gamma'&: \bar V\times \bar U \to Z.
\end{align*}

Fix the parameters $u_1=x_1(1)$ and $u_6=x_6(1)$ in $U_1$ and $U_6$
respectively, {their images $\bar u_1$ and $\bar u_6$ in
$\bar U$,} 
and consider the linear maps
\[\lambda_i:\bar V\to Z\]
given by $\gamma'(x, \bar u_i)$ for $i=1$ or $6$.
The image of $\lambda_1$ is $U_3$ and the kernel is $\bar U_2$ {(use [1,5])}. The image of $\lambda_2$
is $U_4$ and kernel is $\bar U_5$ {(use [2,6])}. So $\bar U_2,U_3,U_4,\bar U_5$ are definable.
Furthermore $\lambda_1$ gives an isomorphism $\lambda_{53}$ of $\bar U_5$ with $U_3$, 
and $\lambda_2$ gives an isomorphism $\lambda_{24}$ of $\bar U_2$ with $U_4$, given  by
\begin{align*}
\lambda_{53}(\bar x_5(a))&=x_3(a);&
\lambda_{24}(\bar x_2(t))=x_4(t).
\end{align*}
We treat these as canonical definable identifications of $\bar U_2$ with $U_4$, and of
$\bar U_5$ with $U_3$.

Now we consider the linear maps $\pi_i:\bar U\to \bar V$ for $i=1$ or $6$ defined by
$\gamma(\bar u_1,x)$ and $\gamma(x, \bar u_6)$. The kernel of $\pi_i$ is
$\bar U_i$ {(by [1,6])} so these groups
are definable. 
Since $\bar U$ and $\bar V$ split definably as $\bar U_1\times \bar U_6$ and $\bar U_2\times \bar U_5$,
the map $\gamma$ induces  maps
\[\mu_j:\bar U_1\times \bar U_6\to \bar U_j\]
with $j=2$ or $5$, given  by
\begin{align*}
\mu_2(\bar x_1(a),\bar x_6(t))&=\bar x_2(-t a^3);&
\mu_5(\bar x_1(a),\bar x_6(t))&=\bar x_5(-t a).
\end{align*}

So with {$a=-1$} this gives a definable isomorphism of $\bar U_6$ and $\bar U_2$ respecting
coordinates, as well as an embedding of $\bar U_6$ into $\bar U_5$ respecting coordinates.
Furthermore with {$t=-1$} we get a definable isomorphism of $\bar U_1$ with $\bar U_5$ 
respecting coordinates, and a definable embedding of $\bar U_1$ into $\bar U_2$ corresponding to the cubing map.

To sum up: we have isomorphic copies of the root groups $U_i$ obtained as definable subquotients of $U$. We have definable identifications of those with even index, respecting coordinates, and of those with odd index, also respecting coordinates. So it is not excessive now to  change notation
and to refer to those with even index as $k$ and to those with odd index as $K$, and to take
the map from $k$ to $K$ as an inclusion, letting $\xi:K\to k$ be the map in the reverse direction, which in coordinates would correspond to cubing.

Modulo these identifications we can simplify the notation (while
recalling that the multiplication is not yet defined in $U$) and
write $m_i(a,t)$ for $\mu_i(\bar x_1(a),\bar x_6(t))$ ($a\in K$, $t\in k$),
and obtain

\[m_2(a, t)=-a^3t, \qquad m_5(a, t)=-at, \qquad  \xi(a)=a^3.\]
with $a\in K$ and $t\in k$. From these maps we can define multiplication on $K$ 
by the condition
$m_2(a,\xi(b))=-\xi(ab)$. We now have $K, k$ with the field structure on $K$ and an inclusion map from $k$ to $K$.

This concludes the proof.
\end{proof}
 
\begin{cor}
\label{Cor:G2:U:biinterpretabiity}
Let 
$(K, k)$ be a pair of fields in characteristic $3$ 
with
\[K^3\le k\le K.\]
Let $U$ be the group $U(k,K)$ in the sense of $\Gtwo(k,K)$ and let $U_i$ be the root groups occurring in $U$,
with the usual numbering.
Then the following structures are bi-interpretable:
\begin{enumerate}
\item $(K,k)$;
\item $(U;U_1,U_6, {u_1,u_6})$;
\item $(U;(U_i)_{1\le i\le 6}, u_1,u_6)$.
\end{enumerate}
\end{cor}

\begin{proof}
This means that we can not only interpret each in the other, but if we apply both interpretations to recover a copy of one of the structures in itself, the result is definably isomorphic to the original structure. As usual, if we start with $(K,k)$ we get an explicit coordinatization of the third structure,
and then we can simply trace the interpretation in coordinates. 

The other direction is a bit less clear since we do in fact require the specified root groups, and not just $U$. The point is that as a set we can identify $U$ definably with the cartesian product
$U_1\times \cdots \times U_6$, but if the root groups were not given to
us this would not be possible.  
\end{proof}

As a corollary we obtain 
\begin{thm}
\label{thm:hexagon}  The group $\Gtwo(k, K)$ is stable
(model-theoretically simple, NTP$_2$, NSOP$_1$,
  \dots) if and only if the pair of fields $(K, k)$ is stable
  (resp., model-theoretically simple, \dots).
\end{thm}

\section{Moufang polygons of indifferent type}
\label{Sec:Indifferent}

Now we come to the case of particular interest, associated with Moufang polygons of indifferent type, or, alternatively, mixing the terminology of Tits and Timmesfeld, the groups
$\PSp_4(L_0,K_0)$  as in Definition
\ref{Def:PSp4(L0,K0)}.

There are some complications at the outset at this level of generality. 
The rank $1$ subgroups of this group have the form $\SL_2(L_0,K_0)$, sitting inside $\SL_2(K)$, but there is a subtlety. One expects to see either $\SL_2(K)$ or $\PSL_2(K)$ here, depending on the root chosen,  but we are in characteristic $2$ so we are free to call the group $\SL_2(K)$, and then recognize the subgroup generated by either $L_0$ or $K_0$ as one of Timmesfeld's groups.
Hence we don't really need to concern ourselves with this point.

A more substantial concern is the nature of the torus in $\PSp_4(L_0,K_0)$, and more generally,
the possible need to enlarge that torus to a larger group $T$ so that $T\PSp_4(L_0,K_0)$ has better definability properties. 
The torus in $\PSp_4(L_0,K_0)$ is generated by rank $1$ root tori coming from $\SL_2(K_0)$ and
$\SL_2(L_0)$, in each case parametrized by the subgroup of $K$ generated by the nonzero elements of $K_0$ or $L_0$, 
respectively (or by the values under the root map giving the multiplicative action on two root groups, for a pair of simple roots).

We elaborate briefly. It is clear that these rank $1$ tori appear in
$\PSp_4(L_0,K_0)$. It is also clear that $\PSp_4$ is generated by the
two rank 1 subgroups associated to simple roots, since they generate $U$
and $U\op$. It follows that if we let $T$ be the torus generated by the
two rank $1$ tori and write out {the} Bruhat decomposition using $TU$ as the group $B$, that the resulting union of double cosets is in fact the full group, and in particular the torus of $\PSp_4(L_0,K_0)$ lies in $B$.

The question that needs to be addressed is, how much of the full torus $T(K)$ of $\PSp_4(K)$ normalizes the group $\PSp_4(L_0,K_0)$.

\begin{lem}\label{Lem:largetorus}
Given a weak indifferent set $(K,L_0,K_0)$, the subgroup of the
full torus $T(K)$ in $\PSp_4(K)$ which 
normalizes $\PSp_4(L_0,K_0)$ is generated by 
$H_\alpha(K)$ and $H_\beta(K_1)$, where $K_1$ is the subfield of $K$ which stabilizes $K_0$.
Here $\alpha,\beta$ are simple roots with $\alpha$ short, and $H_\alpha$, $H_\beta$
are the corresponding rank 1 tori in $T(K)$.
\end{lem}

\begin{proof}
We compute via \cite[Lemma 19]{St}.

The Cartan matrix of this root system is
\[
\begin{bmatrix}
2&-1\\
-2&2
\end{bmatrix}
\]
where the first root is short and the second is long.

This means that the short root torus {$H_\alpha(K)$} operates by squaring (or its inverse) on both root subgroups (for simple roots), and the long root torus $H_\beta(K)$ operates by squaring on $U_\beta(K)$ and 
multiplication (after inversion) on the other root subgroup.

As $K^2\le L_0\le K_0$, the only constraint on
the torus is that the the element $h_\beta(t)$ will occur only if $tK_0=K_0$.
\end{proof}

Notice that $K_1$ contains the field generated by $L_0$, which we know a priori given the structure of the torus in $\PSp_4(L_0,K_0)$---but it is reassuring that it also follows from the definition of
$K_1$ and the hypotheses.

We mentioned earlier that $\PSp_4(L_0,K_0)$ ({see Definition
\ref{Def:PSp4(L0,K0)} and Lemma \ref{Lem:simpleG2Sp4}}) is simple. For definability purposes we may want to consider the group $T\PSp_4(L_0,K_0)$ for some torus normalizing $\PSp_4(L_0,K_0)$ with
better (simpler) definability properties. The goal as always is to obtain simple groups with stable theory of the specified type, or, failing that, at least to find a group $T\PSp_4(L_0,K_0)$ with stable theory, and be content to have a simple commutator subgroup.

Accordingly, it remains only to discuss definability issues in groups of the form
$T\PSp_4(L_0,K_0)$.
We can more or less combine Theorems \ref{Thm:G(k,K)}  and 
\ref{Thm:Rank1:Definability} to get a definability theorem for $T\PSp_4(L_0,K_0)$,
namely Theorem \ref{Thm:Indifferent:Definability:U}
 below. But since we also intend to look more narrowly at definability in $U(k,K)$, we begin with that.

\begin{thm}
  \label{Thm:C2:U(k,K)}
Let $(K; L_0, K_0)$ be a weak indifferent set 
and let $U$ be the group $U(k,K)$ in the sense of $\PSp_4(k,K)$.
Then each of $U$ and $(K_0,L_0,+,*)$ is definable in the other, where
\[a*b=a^2b\]
on $K_0$.
\end{thm}

\begin{proof}
To go from the coordinate system to the group, recall that we can find
a definable field $\tilde K$ so that 
\[({\tilde K^2}\le L_0\le K_0\le \tilde K).\]

As usual $U$ is constructed explicitly from root groups $U_1, U_2, U_3, U_4$
which are 
 copies of $L_0$ or $K_0$, 
 most of which are taken to commute, 
 while
 $U_1, U_4$ satisfy
 \[\tag{[1,4]} [x_1(t),x_4(a)]=x_2(t^2 a)x_3(ta),\]
 where $t$ runs over $K_0$ and $a$ runs over $L_0$, and we use the vector space structure
 on $L_0$ and $K_0$.
 We illustrate the structure of this formula as follows, for reference.
 This is clearly definable in our coordinate system, if we include  $\tilde K$ (to get the action of $L_0$ on $K_0$).

In the converse direction, there is more work to be done.

Let $Z=Z(U)$, {a definable subgroup}. From the commutator formula and the fact that all other pairs of root groups
in $U$ commute, we find that $Z=U_2U_3\simeq U_2\times U_3$.

\smallskip
\noindent \textbf{Claim 1.}  
The subset $U_2 \union U_3$ of $U$ is definable.

\smallskip
We may describe $U_2\union U_3$ as the set of elements $x_2(a)x_3(t)$ with
$a\in L_0$, $t\in K_0$ and at least one coordinate equal to $0$. For $t\in K_0$ and $a\in L_0$ both nonzero, we have
\[[x_1(at^{-1}), x_4(t^2a^{-1}])]=x_2(a)x_3(t).\]
Thus elements of $Z$ which are not in $U_2\union U_3$ {are therefore
  commutators}. On the other hand, by inspection of the commutator
formula [1,4], we see that
the only elements of $U_2\union U_3$ which are {commutators} are those
of the form $[x_1(t),x_4(a)]$ with at least one of $t$ or $a$ equal to $0$; such an element must 
be the identity. This proves the claim. 

\medskip

We will now fix some additional parameters $u_i=x_i(1)\in U_i$ for $i=1,2,3,4$.
We let $\bar U_1$, $\bar U_4$ denote the images of $U_1$, $U_4$ in $\bar
U=U/Z$, and  note that {$\bar U=\bar U_1\times \bar U_4$}.

\smallskip
\noindent \textbf{Claim 2.}  
Relative to the specified parameters, the groups $\bar U_1, U_2, U_3, \bar U_4$ are definable.

By the commutator formula, the centralizer in $U$ of $u_1$ is $U_1Z$, and the image of this group
in $\bar U$ is $\bar U_1$. So $\bar U_1$ is definable, and similarly
$\bar U_4$ is definable.

\smallskip
On the other hand, $U_2$ is the intersection of $u_2(U_2\union U_3)$
with $U_2\union U_3$ 
and so  the group $U_2$ is definable. The group $U_3$ is definable similarly.

This proves the claim.

\smallskip
Now the commutator induces a bilinear map
\[\bar U\times \bar U\to Z\]
and in view of Claim 2, we can interpret the commutation law 
as giving two definable functions of two variables from $\bar U_1\times \bar U_4$ to $U_2$ and $U_3$
respectively. 
These two functions are the maps
\begin{align*}
\mu_2:(\bar x_1(t),\bar x_4(a))\mapsto x_2(t^2a)\;{\in U_2}, &\quad\mbox{(here $t^2a \in L_0$);} \\
\mu_3:(\bar x_1(t),\bar x_4(a))\mapsto x_3(ta)\;{\in U_3}  &\quad \mbox{(here $ta \in K_0$).}
\end{align*}

With $t=1$ the map $\mu_2$ gives an isomorphism $\bar U_4\to U_2$ respecting coordinates.

With $a=1$  the map $\mu_3$ gives an isomorphism $\bar U_1\to U_3$ respecting 
coordinates. 

With these definable isomorphisms available, we identify $\bar U_1$ with
$U_3$ and {with} the additive group {of}
$K_0$. Likewise, we identify $\bar U_4$ with $U_2$ and  the additive
group of 
$L_0$. Now $K_0$, $L_0$ are interpreted in $U$ and the maps $\mu_2$, $\mu_3$ become
definable maps $m_2,m_3$ from $K_0\times L_0$ to $L_0$ and $K_0$ respectively, 
satisfying the following conditions:
\begin{align*}
m_2(t, a)&=t^2a;\\
m_3(t, a)&=ta.
\end{align*}
With $t=1$ the map $m_3$ defines the embedding of $L_0$ into $K_0$, so we now view
$L_0$ as a subset of $K_0$, and both maps $m_2$ and $m_3$ now have image
in $K_0$ {(identified with $U_3$)}. 
In other words,
the definable structure present is
\[(K_0; L_0,+, m_2,m_3).\]
With $a=1$, $m_2$ defines the squaring map on $K_0$.
Now for $a, b\in K_0$, we have $a^2\in L_0$, and {$m_3(b,
a^2)=a^2b$}. 
\end{proof}

\begin{cor}
Let 
Let $(K; L_0, K_0)$ be a weak indifferent set 
and let $U$ be the group $U(k,K)$ in the sense of $\PSp_4(k,K)$, with root subgroups $U_i$ ($1\le i\le 4$)
numbered as usual.
Then the following structures are bi-interpretable:
\begin{enumerate}
\item $(K_0,L_0;+,*)$ where $a*b=a^2 b$ for $a,b\in K_0$.
\item $(U;+, (U_i)_{1\le i\le 4}, u_1,\ldots,u_4)$ with $U_i$ the root
  groups and $u_i=x_i(1)\in U_i$.
\end{enumerate}
\end{cor}

\begin{proof}
As in the proof of 
Corollary \ref{Cor:G2:U:biinterpretabiity}.
\end{proof}

\begin{thm}
\label{Thm:Indifferent:Definability:U}
Let $(K; L_0, K_0)$ be a weak indifferent set, $T(K)$ a maximal torus of $\PSp_4(K)$,
and $T$ a subgroup of $T(K)$ normalizing the group
$\PSp_4(L_0,K_0)$
and containing $T(K)\intersect \PSp_4(L_0,K_0)$.

Let $\calm$ be the structure
\[(K_0;L_0,T,+,\mu)\]
consisting of the group $K_0$ with the subset $L_0$, the abstract group $T$ with its multiplication,
and the following additional structure:
\begin{enumerate}
\item  the map $\mu:K_0\times K_0\to K_0$ defined by $\mu(a,b)=a^2b$;
\item actions of $T$ on $K_0$ and on $L_0$ which
correspond to the actions of $T$ on two root subgroups $U_\alpha$, $U_\beta$
with $\alpha,\beta$ the two simple roots, where $\alpha$ is short and $\beta$ is long.
\end{enumerate}

Then the group $G=T\PSp_4(L_0,K_0)$ is interdefinable with $\calm$.

In particular, $G$ is stable if and only if $\calm$ is stable. 
\end{thm}

Note that by Lemma \ref{Lem:DefinableFields} we may also include fields $\tilde K, \tilde L$ with
\[\tilde K^2\le L_0\le \tilde L\le K_0\le \tilde K,\]
where $\tilde K$ contains $K$, but to do so we need to choose fields which are definable in the 
given structure. 

Typically $T$ is a product of ``root tori'' (intersections with root tori in $T(K)$ and this can be simplified further to just take the action of a given root torus on the corresponding root group. 
But we prefer to give the statement in its general form.

\begin{proof}

We repeat the proof of Theorem \ref{Thm:G(k,K)} 
with minor adjustments. 

1. First we show that $G$ is definable from the coordinatizing structure $\calm$.

We will take a definable field $\tilde K$ with 
\[\tilde K^2\le L_0\le K_0\le K\le \tilde K\]
and work in the larger group $G(\tilde K)$ as the ``ambient'' group.

We build the group $U(k,K)$ as a subgroup of $\PSp_4(\tilde K)$ using the coordinate system 
$(\tilde K;L_0,K_0)$. We use the action of $T$ on $L_0$ and $K_0$ to define an action of $T$ on the root groups $U_\alpha(K_0)$ and $U_\beta(L_0)$ for $\alpha,\beta$ simple,
$\alpha$ short, $\beta$ long, which determines the action on $U$ by the commutator formula; that is,
the two given root groups generate $U(k,K)$ in a bounded way. One can then find a unique embedding of $T$ into $T(\tilde K)$ which respects that action. So now $B$ is a definable subset (in an appropriate language) in $\PSp_4(\tilde K)$. The Bruhat decomposition therefore gives $G$ as a definable subset of $\PSp_4(\tilde K)$, and the multiplication law is inherited, 
so the group $G$ is definable.\\

2. For the converse, we first find that $U$ and $T$  are definable using
the definability of centralizers---here is the argument. Working in
$\Gg(K^{alg})$ for $\Gg$ an
algebraic group, we know
that the descending chain condition on centralizers holds, i.e., that
given $X\subset \Gg(K^{alg})$, there is some finite $X_0\subset X$ such
that $C_{\Gg(K^{alg})}(X)=C_{\Gg(K^{alg})}(X_0)$. But this property is
  preserved by going to subgroups: i.e., if $X\subset \PSp_4(k,K)$, then
  $C_{\PSp_4(k,K)}(X)$ is definable. In particular the torus $T_0(k,K)$
  is self-centralizing in $\PSp_4(k,K)$, hence definable, and the centralizer $Z$ of
  $U(k,K)$ is the center of $U(k,K)$. So $T_0(k,K)$ and $Z$ are
  definable. 

\smallskip\noindent
\textbf{Claim}. 
Each root subgroup  $U_\alpha$ is equal to $CC(U_\alpha)$,
hence definable.

\smallskip\noindent
{\em Proof of the claim}.  
First, we can suppose the root subgroup is in the center of $U$,
{because root groups of the same length are conjugate (see
Appendix I.15 in  \cite{St}, or Lemma C, section 10.4 in \cite{Hum-Lie} in the classical case)}. Then its centralizer contains $U$,
and contains the rank $1$ component 
{$L_\beta=\gen{U_\beta,U_{-\beta}}$}, with $\pm\beta$ the unique root orthogonal to 
$\alpha$ {so that $L_\beta\leq C_G(U_\alpha)$}. Hence the double
  centralizer of $U_\alpha$ is contained in the center of $U$ and
  commutes with $L_\beta$. Looking at Lemma~19(c) in \cite{St}, we see
  that the action of a typical element  $h_\beta(t)$ of the torus of
  $L_\beta$ acts on a root subgroup $U_\gamma$ as multiplication by
  $t^{\langle \beta,\gamma\rangle}$, where $\langle \beta,\gamma\rangle$
  denotes the inner product of the roots
  $\beta,\gamma$ in $\mathbb{R}^2$. In our case $\alpha$ is the only positive root
  which is orthogonal to $\beta$, which proves our assertion.
 \qed

\smallskip
If $L_\alpha=\gen{U_{\pm \alpha}}$, we can
also define the group $TL_\alpha$. Indeed, we have the group $B_\alpha=TU_\alpha$
since we have $U_\alpha$, and then $TL_\alpha$ is definable via its Bruhat decomposition.
If $k=L_0$ or $K_0$ is the parametrizing group for $U_\alpha$, this group gives us
$(k,T)$ with $T$ represented via its action on $k$. 

As we have the group $U$ we have the remaining structure on the coordinate system by 
Theorem \ref{Thm:Indifferent:Definability:U}.
\end{proof}

\section{Concluding Remarks}
While the paper is certainly not self-contained, we have developed some
parts of the group theoretic analysis in large detail, as it casts
considerable light on the definability issues and other important
matters. Much of this is of the sort found in Steinberg's notes
\cite{St} or, from a radically different point of view, in the book of
Tits and Weiss \cite{TW}. 

\medskip
Once again, Theorem~\ref{thm2} (as well as Theorem \ref{thm3}) applies to give a number of stable simple groups acting as automorphism groups of Moufang polygons of indifferent type, or in some cases stable groups with simple commutator subgroup acting as autmorphism groups of the same kind. In particular, in some cases adding a nondefinable field to the coordinatizing structure will provide additional examples which are in fact simple.

\end{document}